%% file: Arxiv_submission.tex
\documentclass[11pt]{article}

\usepackage{amsmath}
\usepackage{mathrsfs}
\usepackage{amssymb}
\usepackage{subcaption}
\usepackage{graphicx}
\usepackage{algorithm}
\usepackage{algorithmic}
\usepackage{natbib}
\usepackage{lipsum}
\usepackage{makecell}
\usepackage{enumitem}
\usepackage{listings}
\usepackage[most]{tcolorbox}
\usepackage{multirow}
\usepackage{backref}
\usepackage{wrapfig}

\usepackage{natbib}

\usepackage{url}
\usepackage[colorlinks,linkcolor=magenta,citecolor=blue, pagebackref=true,backref=true]{hyperref}
\renewcommand*{\backrefalt}[4]{%
    \ifcase #1 \footnotesize{(Not cited.)}%
    \or        \footnotesize{(Cited on page~#2.)}%
    \else      \footnotesize{(Cited on pages~#2.)}%
    \fi}

\usepackage{amsthm}
\usepackage{amsfonts}
\usepackage{graphicx}
\usepackage{epstopdf}
\usepackage{tikz}
\usepackage{mathtools}

\usepackage[capitalize]{cleveref}
\Crefname{lemma}{Lemma}{Lemmas}
\Crefname{remark}{Remark}{Remarks}
\Crefname{figure}{Figure}{Figures}
\Crefname{enumi}{Lemma}{Lemmas}
    
\newtheorem{theorem}{Theorem}[section]
\newtheorem{corollary}[theorem]{Corollary}
\newtheorem{lemma}[theorem]{Lemma}

\newtheorem{definition}{Definition}[section]
\newtheorem{example}{Example}[section]

\newtheorem{remark}[theorem]{Remark}
\newtheorem{assumption}[theorem]{Assumption}

\newcommand{\BB}{\mathbb{B}}
\newcommand{\EE}{\mathbb{E}}

\newcommand{\Exp}{\textnormal{Exp}}

\newcommand{\grad}{\textnormal{grad}}

\newcommand{\argmin}{\mathop{\rm argmin}}

\newcommand{\MCal}{\mathcal{M}}
\newcommand{\NCal}{\mathcal{N}}

\newcommand{\SCal}{\mathcal{S}}

\newcommand{\br}{\mathbb{R}}

\newcommand{\ba}{\begin{array}}
\newcommand{\ea}{\end{array}}

\begin{document}


\begin{center}

{\bf{\LARGE{Last-Iterate Convergence of Adaptive Riemannian \\ [\medskipamount] Gradient Descent for Equilibrium Computation}}}

\vspace*{.2in}
{\large{ \begin{tabular}{c}
Yang Cai$^\dagger$ \and Michael I. Jordan$^\circ$ \and Tianyi Lin$^\ddagger$ \and Argyris Oikonomou \and Emmanouil V. Vlatakis-Gkaragkounis$^\diamond$ \\
\end{tabular}
}}

\vspace*{.2in}

\begin{tabular}{c}
Department of Computer Science, Yale University$^\dagger$ \\
Department of Electrical Engineering and Computer Sciences, UC Berkeley$^\circ$ \\
Department of Industrial Engineering and Operations Research, Columbia University$^\ddagger$ \\ 
Department of Computer Sciences, University of Wisconsin-Madison$^\diamond$
\end{tabular}

\vspace*{.2in}

\today\\ [.2cm]

\vspace*{.2in}

\begin{abstract}
Equilibrium computation on Riemannian manifolds provides a unifying framework for numerous problems in machine learning and data analytics. One of the simplest yet most fundamental methods is Riemannian gradient descent (RGD). While its Euclidean counterpart has been extensively studied, it remains unclear how the manifold curvature affects RGD in game-theoretic settings. This paper addresses this gap by establishing new convergence results for \textit{geodesic strongly monotone} games. Our key result shows that RGD attains last-iterate linear convergence in a \textit{geometry-agnostic} fashion, a key property for applications in machine learning. We extend this guarantee to stochastic and adaptive variants---SRGD and FARGD---and establish that: (i) the sample complexity of SRGD is geometry-agnostic and optimal with respect to noise; (ii) FARGD matches the convergence rate of its non-adaptive counterpart up to constant factors, while avoiding reliance on the condition number. Overall, this paper presents the first geometry-agnostic last-iterate convergence analysis for games beyond the Euclidean settings, underscoring the surprising power of RGD---despite its simplicity---in solving a wide spectrum of machine learning problems.
\end{abstract}

\end{center}

\input{sec/intro}

\input{sec/prelim}

\input{sec/game}

\input{sec/algorithm}

\input{sec/exp}

\input{sec/conclu}


\bibliographystyle{plainnat}
\bibliography{ref}

\newpage
\appendix
\input{sec/app}

\end{document}

%% file: sec/intro.tex
\section{Introduction}
Game-theoretic concepts underpin a broad spectrum of machine learning (ML) domains, ranging from single-agent optimization~\citep{Beck-2003-Mirror, Nemirovski-2004-Prox, Petcu-2005-Scalable, Lan-2012-Validation} and feature selection~\citep{Dikkala-2020-Minimax, Jimenez-2021-Novel} to contemporary applications such as adversarial training~\citep{Goodfellow-2014-Generative, Kumar-2017-Semi, Madry-2018-Towards}, multi-agent dynamics~\citep{Bucsoniu-2010-Multi, Vlatakis-2020-Noregret}, and AI alignment~\citep{Munos-2024-Nash,Swamy-2024-Minimaximalist,Liu-2024-Comal}. A prevailing methodological paradigm is to characterize the task's optimal solution as a Nash equilibrium~\citep{Nash-1951-noncooperative, Rosen-1965-Existence} of a suitably defined game. Yet, closer inspection of the literature reveals that the convergence guarantees of many state-of-the-art algorithms rely on implicit structural assumptions---most notably, that the underlying game is ``nearly'' \textit{concave} or \textit{monotone} and that the feasible strategy set is \textit{convex}~\citep{Facchinei-2007-Finite, Diakonikolas-2021-Efficient, Cai-2022-Finite, Pethick-2023-Solving}.

These structures are, to some extent, introduced to serve the needs of theoretical analysis. In various real-world applications, ranging from optimal transport~\citep{Lin-2020-Projection, Huang-2021-Riemannian} and environmental engineering~\citep{Oliger-1978-Theoretical, Kinderlehrer-2000-Introduction} to robotics and multi-agent systems involving complex interactions~\citep{Spong-1987-Integral, Quinlan-1994-Efficient, Berenson-2009-Manipulation, Bhattacharya-2013-Distributed}, the games often exhibit nonconvexity, either in the utility functions or in the strategy sets. Such nonconvexities pose significant challenges for conventional approaches, as highlighted in~\citet{Mertikopoulos-2018-Cycles, Vlatakis-2019-Poincare, Daskalakis-2021-Complexity}, and frequently necessitate the use of ad hoc heuristics drawn from optimization literature~\citep{Toscano-2010-New, Diamond-2018-General}. This raises a key research challenge: to identify geometric structures of utility functions or strategy sets that extend beyond near-concavity and monotonicity, and to design algorithms that systematically leverage such structure. Recent research has sought to establish geometric foundations for utility functions~\citep{Mertikopoulos-2019-Learning, Mazumdar-2020-Gradient, Hsieh-2021-Adaptive, Vlatakis-2021-Solving, Sakos-2023-Exploiting, Mertikopoulos-2024-Unified, Mazumdar-2025-Finding}. Complementing this direction, another line of work has recast constraints in game theory and equilibrium computation within the framework of \emph{Riemannian manifolds}~\citep{Jordan-2022-First, Han-2023-Riemannian, Zhang-2023-Sion, Hu-2024-Extragradient, Li-2025-Riemannian}.

This geometric framework, which commonly employs Hadamard and Stiefel manifolds, has been applied to a wide range of statistical tasks, including online principal component analysis~\citep{Lee-2022-Fast}, diffusion tensor data processing~\citep{Lenglet-2006-Geometric}, and maximum likelihood estimation for heavy-tailed, non-Gaussian distributions~\citep{Wiesel-2012-Geodesic}. Extending this framework to game-theoretic domains, however, remains challenging, particularly in settings involving nonconvex strategy sets. Examples include robotic systems, where the feasible joint angles of robotic arms lie on the SO(3) rotation manifold~\citep{Hong-2020-Real, Watterson-2020-Trajectory}, and environmental pollution control games, where the physical transmission of pollutants across bounded surface areas yields manifold constraints~\citep{Nagurney-1997-General, Scrimali-2010-Variational}.

These domains present two fundamental challenges: \textit{the existence of equilibria} and \textit{the convergence of game dynamics}. Unlike optimization problems, which admit a solution over compact domains, the existence of Nash equilibria often hinges on fixed-point theorems~\citep{Brouwer-1911-Abbildung, Kakutani-1941-Generalization}, which require convexity of the strategy sets. Formally, the goal is to identify a joint action $x^\star = (x_1^\star, \dots, x_N^\star)$ such that
\begin{equation*}
u_i(x^\star) \geq u_i(x_i; x_{-i}^\star), \quad \textnormal{for all } x_i \in \MCal_i \textnormal{ and all } i,
\end{equation*}
where $(x_i; x_{-i})$ denotes the joint action of the $i$-th player and all other players, $\{\MCal_i\}_{1 \leq i \leq N}$ denotes nonconvex strategy sets, and $\{u_i\}_{1 \leq i \leq N}$ denote utility functions. Existence results for Nash equilibria in geodesically concave games over Hadamard manifolds were first obtained in~\citet{Nemeth-2003-Variational} and extended in~\citet{Li-2009-Existence} through the framework of variational inequalities. Recently,~\citet{Zhang-2023-Sion} have proved a strong duality theorem for minimax Riemannian optimization applicable to a broad class of Riemannian manifolds with \textit{unique geodesics}.

A general understanding of the convergence properties of game dynamics on manifolds remains limited. To bridge this gap, a line of research seeks to generalize classical algorithms, including gradient descent ascent (GDA) and extragradient (EG) methods~\citep{Jin-2022-Understanding, Jordan-2022-First, Huang-2023-Gradient, Zhang-2023-Sion}. Existing works focus on min-max formulations. Notably,~\citet{Jordan-2022-First} has shown that the Riemannian corrected extra-gradient (RCEG) method achieve last-iterate linear convergence for solving geodesically strongly convex-concave optimization problems, thus extending the convergence guarantees known in the Euclidean setting to the Riemannian regime.

Proceeding to the general setting of multi-player games where strategy spaces lie on manifolds, several fundamental challenges emerge. First, we need to clarify the appropriate solution concept in this inherently nonconvex games. Second, it is essential to identify structures of these nonconvex games that render the appropriate solution concepts computationally tractable.\footnote{Without any structural assumptions, the problem is computationally intractable~\citep{Daskalakis-2021-Complexity}.} Finally, it is desirable to design game dynamics that are uncoupled and, more importantly, operate in a \textit{geometry-agnostic} fashion, indicating the independence of any \textit{sectional curvature} parameters. 

\paragraph{Related works} Research on min--max optimization over Riemannian manifolds has largely centered around two fundamental questions: (i) the existence and uniqueness of equilibria, and (ii) the efficient computation of such equilibria. The first line of results was established for Hadamard manifolds~\citep{Komiya-1988-Elementary, Kristaly-2014-Nash, Park-2019-Riemannian} and later generalized to finite-dimensional manifolds with unique geodesics~\citep{Zhang-2023-Sion}. Algorithms for equilibrium computation were developed for hyperbolic Hadamard manifolds with negative curvature~\citep{Ferreira-2005-Singularities, Li-2009-Monotone, Wang-2010-Monotone}, though these methods only guaranteed asymptotic convergence. Non-asymptotic convergence analysis has appeared only recently for gradient-based algorithms in Riemannian min-max optimization. In geodesically smooth and strongly convex--strongly concave settings, the best known method achieves last-iterate and \textit{geometry-aware} global convergence~\citep{Jordan-2022-First}. In geodesically convex--concave settings, existing algorithms rely on \textit{geometry-aware} step sizes and achieve sublinear convergence rates: $O(1/T)$ in the time-average sense~\citep{Zhang-2023-Sion,Jordan-2022-First} and $O(1/\sqrt{T})$ in the last-iterate sense~\citep{Hu-2024-Extragradient}. To the best of our knowledge, our work provides the first results that bridge this gap in equilibrium computation for games on Riemannian manifolds, demonstrating the possibility of achieving last-iterate and \textit{geometry-agnostic} convergence guarantees for gradient-based algorithms.

Independent and concurrent work of \citet{Martinez-2025-Accelerated} investigates Riemannian min–max optimization on Hadamard manifolds. Their results most relevant to ours establish the existence of equilibria in geodesically quasi-convex-quasi-concave settings under the assumption of unique geodesics. This contrasts with our results, which address equilibrium existence in general strongly monotone games (cf. Assumption~\ref{Assumption:main-1}) without requiring uniqueness of geodesics. In addition, they propose novel double-loop accelerated algorithms tailored to their specific setting, although their convergence guarantees remain geometry-dependent, relying explicitly on the sectional curvature parameters of the manifold.

\paragraph{Contributions} This paper advances the theory of game-theoretical learning on Riemannian manifolds by identifying a class of games that admit simple algorithms with linear, geometry-agnostic, and last-iterate convergence guarantees. Although our framework covers practical applications, the main contributions are theoretical and can be summarized as follows:
\begin{enumerate}
\item We introduce and study a class of geodesically strongly monotone games, where each player chooses an action on a Riemannian manifold and the players’ loss functions jointly satisfy the geodesic monotonicity condition. We further prove that any geodesically strongly monotone game admits a Nash equilibrium, generalizing the equilibrium existence results of~\citet{Zhang-2023-Sion} to multi-player general-sum settings.
\item We analyze the convergence properties of Riemannian Gradient Descent (RGD) for geodesically strongly monotone games. We show that RGD with a constant step size, which is chosen independently of the manifold curvature, can achieve linear, geometry-agnostic, and last-iterate convergence to the Nash equilibrium. We also extend these results to stochastic and fully adaptive variants of RGD.
\end{enumerate}
From a technical viewpoint, our results establish the first linear, geometry-agnostic and last-iterate convergence rates for games on manifolds, covering Riemannian min-max optimization as a special case and eliminating the need for unique geodesic assumptions. This resolves an open question from~\citet{Jordan-2022-First} and opens new directions for Riemannian game dynamics. Our key innovation lies in a convergence analysis that departs from previous methods relying on \emph{Trigonometric Comparison Inequalities} (TCIs), which inherently induce curvature dependence and impose additional constraints on the initialization. Instead, we provide a generalized descent lemma for Riemannian strongly monotone games and extend it to a stochastic setting using Young's inequality. Although this leads to an additional factor that depends on the condition number of the underlying function in the convergence rate, it guarantees that the resulting bounds remain geometry-agnostic. One limitation is that we assume no boundaries for the manifolds. Extending our results to manifolds with boundaries necessitates well-defined projection mappings and lies beyond the scope of this work.

\paragraph{Organization} In Section~\ref{sec:prelims}, we present background materials on manifold geometry. In Section~\ref{sec:games}, we provide the definitions of geodesically strongly monotone games and define our performance function. In Section~\ref{sec:results}, we prove our main results for RGD and extend these results to stochastic and adaptive variants. In Section~\ref{sec:exp}, we present empirical results complementing our theoretical results. We conclude in Section~\ref{sec:conclu}. 

%% file: sec/prelim.tex
\section{Preliminaries}\label{sec:prelims}
An $n$-dimensional manifold $\MCal$ is a topological space where any point has a neighborhood that is homeomorphic to $n$-dimensional Euclidean space. For each $x \in \MCal$, each tangent vector is tangent to all of parametrized curves passing through $x$ and the tangent space $T_x \MCal$ of a manifold $\MCal$ at $x$ is defined as the set of all of tangent vectors. A Riemannian manifold $\MCal$ is a smooth manifold that is endowed with a smooth (``Riemannian") metric $\langle \cdot, \cdot\rangle_x$ on the tangent space $T_x \MCal$ for each point $x \in \MCal$. The inner metric induces a norm $\|\cdot\|_x$ on the tangent spaces.   

\begin{figure}[!t] 
\begin{center} 
\scalebox{0.4}{\input{figs/tikz-a}} 
\end{center}     
\caption{Illustrations of tangent spaces, parallel transport and exponential maps.} \label{fig:my_label} 
\end{figure}

A geodesic can be seen as the generalization of an Euclidean linear segment and is modeled as a smooth curve, $\gamma: [0, 1] \mapsto \MCal$, which is locally a distance minimizer. In contrast to the Euclidean setting, $x$ and $v = \grad f(x)$ do not lie in the same space, since $\MCal$ and $T_x \MCal$ respectively are distinct entities. The interplay between these dual spaces typically is carried out via the \textit{exponential maps}. An exponential map at a point $x \in \MCal$ is a mapping from the tangent space $T_x \MCal$ to $\MCal$. In particular, $y=\Exp_x(v) \in \MCal$ is defined such that there exists a geodesic $\gamma$ satisfying $\gamma(0) = x$, $\gamma(1) = y$ and $\dot{\gamma}(0) = v$. 

In contrast again to Euclidean spaces, we cannot compare the tangent vectors at different points $x, y \in \MCal$ since these vectors lie in different tangent spaces. To resolve this issue, it suffices to define a transport mapping that moves a tangent vector along the geodesics and also preserves the length and Riemannian metric $\langle \cdot, \cdot\rangle_x$; indeed, we can define a parallel transport $\Gamma_x^y: T_x \MCal \mapsto T_y \MCal$ such that the inner product between any $u, v \in T_x \MCal$ is preserved; i.e., $\langle u, v\rangle_x = \langle \Gamma_x^y(u), \Gamma_x^y(v)\rangle_y$. 

Given a geodesic $\gamma: [0, 1] \mapsto \MCal$ on a Riemannian manifold $\MCal$, the length of the geodesic $L(\gamma)$ is obtained by integrating the norm of its tangent vector along the curve. In other words, $L(\gamma) = \int_0^1 \|\dot{\gamma}(t)\|_{\gamma(t)} dt$ where $\dot{\gamma}(t)$ is the tangent vector of $\gamma$ at the point $\gamma(t)$. The distance between $x$ and $y$ is $d_\MCal(x, y) = \inf_{\gamma: \gamma(0)=x, \gamma(1)=y} L(\gamma)$. For point $x \in \MCal$ and $D > 0$, we define $\BB_\MCal(x; D)= \{x' \in \MCal: d_\MCal(x, x') \leq D\}$ as the set of points with distance at most $D$ from $x$.

For $x, y \in \MCal$, a geodesic $\gamma: [0,1] \mapsto \MCal$ with $\gamma(0)=x$, $\gamma(1)=y$, and $\grad f(x) \in T_x \MCal$ is a Riemannian gradient of $f$ at a point $x$, the following statements hold true, 
\begin{itemize}
\item $f$ is geodesically $\mu$-strongly convex if $f(y) \geq f(x) + \langle \grad f(x), \dot{\gamma}(0) \rangle_x + \tfrac{\mu}{2}\|\dot{\gamma}(0)\|_x^2$.  
\item $f$ is geodesically $\ell$-smooth if $\|\grad f(x) - \Gamma_y^x \grad f(y)\|  \leq  \ell d_\MCal(x, y)$.  
\end{itemize}
We make use of the following properties of the parallel transport $\Gamma$: 
\begin{enumerate}
\item For $0 \leq t_1 \leq t_2 \leq t_3 \leq 1$ and $u \in T_{\gamma(t_1)}\MCal$, the following statements hold true, 
\begin{equation*}
\Gamma_{\gamma(t_2)}^{\gamma(t_3)} \Gamma_{\gamma(t_1)}^{\gamma(t_2)}(u)= \Gamma_{\gamma(t_1)}^{\gamma(t_3)}(u), \quad \Gamma_{\gamma(t_1)}^{\gamma(t_2)} \Gamma_{\gamma(t_2)}^{\gamma(t_1)}(u)= u.
\end{equation*}
\item For $0 \leq t_1 \leq t_2 \leq 1$, let $\zeta(t) = \gamma(t\cdot (t_2 - t_1) + t_1)$ be the restriction of $\gamma$ on $[t_1,t_2]$. Then, for $v \in T_{\gamma(t_1)} \MCal$, we have $\Gamma_{\gamma(t_1)}^{\gamma(t_2)}(v) =  \Gamma_{\zeta(0)}^{\zeta(1)}(v)$ and $\dot{\zeta}(t) = (t_2-t_1) \cdot \Gamma_{\gamma(t_1)}^{\gamma(t\cdot (t_2 - t_1) + t_1)}(v)$.
\end{enumerate}
We derive an alternative characterization of geodesically strongly convex functions using the above properties of $\Gamma$. For $x, y \in \MCal$, a geodesic $\gamma: [0,1] \mapsto \MCal$ with $\gamma(0)=x$, $\gamma(1)=y$, and with $\grad f(x) \in T_x \MCal$ is a Riemannian gradient of $f$ at a point $x$, a function $f$ is geodesically $\mu$-strongly convex if the following statement holds true, 
\begin{equation*}
f(\gamma(t_2)) \geq f(\gamma(t_1)) + (t_2-t_1) \langle\grad f(\gamma(t_1)), \dot{\gamma}(t_1)\rangle_{\gamma(t_1)} + \tfrac{\mu}{2}(t_2-t_1)^2\|\dot{\gamma}(t_1)\|_{\gamma(t_1)}^2, 
\end{equation*}
for $0 \leq t_1 \leq t_2 \leq 1$. This inequality follows by using the sub-geodesic $\zeta(t) = \gamma(t \cdot (t_2-t_1)+t_1)$ in the original definition of geodesically strongly convex functions and $\dot{\zeta}(0) = (t_2-t_1) \cdot \dot{\gamma}(t_1)$.

%% file: figs/tikz-a.tex
\tikzset{every picture/.style={line width=0.75pt}} 

\begin{tikzpicture}[x=0.75pt,y=0.75pt,yscale=-1,xscale=1]

\draw  [fill={rgb, 255:red, 155; green, 155; blue, 155 }  ,fill opacity=0.26 ] (95.8,179.13) -- (217.35,72.2) -- (197.14,155.36) -- (75.59,262.29) -- cycle ;
\draw    (150.41,165.45) -- (97.26,212.21) ;
\draw [shift={(95.76,213.53)}, rotate = 318.66] [color={rgb, 255:red, 0; green, 0; blue, 0 }  ][line width=0.75]    (10.93,-3.29) .. controls (6.95,-1.4) and (3.31,-0.3) .. (0,0) .. controls (3.31,0.3) and (6.95,1.4) .. (10.93,3.29)   ;
\draw [shift={(150.41,165.45)}, rotate = 138.66] [color={rgb, 255:red, 0; green, 0; blue, 0 }  ][fill={rgb, 255:red, 0; green, 0; blue, 0 }  ][line width=0.75]      (0, 0) circle [x radius= 3.35, y radius= 3.35]   ;
\draw  [fill={rgb, 255:red, 155; green, 155; blue, 155 }  ,fill opacity=0.26 ] (392.86,143.12) -- (476.06,293.97) -- (381.93,247.5) -- (298.73,96.65) -- cycle ;
\draw    (391.72,200.64) -- (355.28,134.56) ;
\draw [shift={(354.31,132.81)}, rotate = 61.12] [color={rgb, 255:red, 0; green, 0; blue, 0 }  ][line width=0.75]    (10.93,-3.29) .. controls (6.95,-1.4) and (3.31,-0.3) .. (0,0) .. controls (3.31,0.3) and (6.95,1.4) .. (10.93,3.29)   ;
\draw [shift={(391.72,200.64)}, rotate = 241.12] [color={rgb, 255:red, 0; green, 0; blue, 0 }  ][fill={rgb, 255:red, 0; green, 0; blue, 0 }  ][line width=0.75]      (0, 0) circle [x radius= 3.35, y radius= 3.35]   ;
\draw    (198,124) .. controls (237.6,94.3) and (283.08,92.04) .. (322.8,121.11) ;
\draw [shift={(324,122)}, rotate = 216.87] [color={rgb, 255:red, 0; green, 0; blue, 0 }  ][line width=0.75]    (10.93,-3.29) .. controls (6.95,-1.4) and (3.31,-0.3) .. (0,0) .. controls (3.31,0.3) and (6.95,1.4) .. (10.93,3.29)   ;
\draw  [draw opacity=0][fill={rgb, 255:red, 208; green, 2; blue, 27 }  ,fill opacity=0.13 ] (109,266.5) .. controls (109,266.5) and (109,266.5) .. (109,266.5) .. controls (109,185.04) and (175.04,119) .. (256.5,119) .. controls (337.1,119) and (402.6,183.65) .. (403.98,263.93) -- (256.5,266.5) -- cycle ; \draw   (109,266.5) .. controls (109,266.5) and (109,266.5) .. (109,266.5) .. controls (109,185.04) and (175.04,119) .. (256.5,119) .. controls (337.1,119) and (402.6,183.65) .. (403.98,263.93) ;  
\draw  [draw opacity=0][fill={rgb, 255:red, 208; green, 2; blue, 27 }  ,fill opacity=0.13 ] (505.07,259.26) .. controls (507.54,182.33) and (571.95,120.7) .. (651.04,120.7) .. controls (731.69,120.7) and (797.08,184.79) .. (797.08,263.85) -- (651.04,263.85) -- cycle ; \draw   (505.07,259.26) .. controls (507.54,182.33) and (571.95,120.7) .. (651.04,120.7) .. controls (731.69,120.7) and (797.08,184.79) .. (797.08,263.85) ;  
\draw  [fill={rgb, 255:red, 155; green, 155; blue, 155 }  ,fill opacity=0.26 ] (737.02,161.85) -- (542.51,161.03) -- (644.08,87.12) -- (838.6,87.93) -- cycle ;
\draw    (651.35,121.7) -- (761.16,122.17) ;
\draw [shift={(763.16,122.18)}, rotate = 180.24] [color={rgb, 255:red, 0; green, 0; blue, 0 }  ][line width=0.75]    (10.93,-3.29) .. controls (6.95,-1.4) and (3.31,-0.3) .. (0,0) .. controls (3.31,0.3) and (6.95,1.4) .. (10.93,3.29)   ;
\draw [shift={(651.35,121.7)}, rotate = 0.24] [color={rgb, 255:red, 0; green, 0; blue, 0 }  ][fill={rgb, 255:red, 0; green, 0; blue, 0 }  ][line width=0.75]      (0, 0) circle [x radius= 3.35, y radius= 3.35]   ;
\draw    (651.35,121.7) .. controls (706.72,131) and (744.56,155.89) .. (746,172) ;
\draw [shift={(746,172)}, rotate = 84.88] [color={rgb, 255:red, 0; green, 0; blue, 0 }  ][fill={rgb, 255:red, 0; green, 0; blue, 0 }  ][line width=0.75]      (0, 0) circle [x radius= 3.35, y radius= 3.35]   ;

\draw (122.66,241) node [anchor=north west][inner sep=0.75pt]  [font=\Large]  {$\MCal$};
\draw (249,55.4) node [anchor=north west][inner sep=0.75pt]  [font=\Large]  {$\Gamma_x^y$};
\draw (657.9,173.71) node [anchor=north west][inner sep=0.75pt]  [font=\Large]  {$y = \Exp_x(v)$};
\draw (640.22,95.07) node [anchor=north west][inner sep=0.75pt]  [font=\Large]  {$x$};
\draw (574.36,135.31) node [anchor=north west][inner sep=0.75pt]  [font=\Large]  {$T_x \MCal$};
\draw (730.51,103.0) node [anchor=north west][inner sep=0.75pt]  [font=\Large]  {$v$};
\draw (772.66,240) node [anchor=north west][inner sep=0.75pt]  [font=\Large]  {$\MCal$};
\draw (390.3,146.13) node [anchor=north west][inner sep=0.75pt]  [font=\Large,rotate=-70.9]  {$v$};
\draw (392.09,205.13) node [anchor=north west][inner sep=0.75pt]  [font=\Large,rotate=-67.41]  {$y$};
\draw (440.72,224.67) node [anchor=north west][inner sep=0.75pt]  [font=\Large,rotate=-67.41]  {$T_y \MCal$};
\draw (153.72,157.41) node [anchor=north west][inner sep=0.75pt]  [font=\Large,rotate=-319.09]  {$x$};
\draw (161.01,120.05) node [anchor=north west][inner sep=0.75pt]  [font=\Large,rotate=-319.09]  {$T_x \MCal$};
\draw (99.4,180.74) node [anchor=north west][inner sep=0.75pt]  [font=\Large,rotate=-324.55]  {$u$};
\end{tikzpicture}

%% file: sec/game.tex
\section{Geodesically Monotone Games}\label{sec:games}
We consider games where the action of player $i$ lies on a Riemannian manifold $\MCal_i$, and extend monotone games to this more general setting. 

\begin{definition}\label{def:game}
A game on Riemannian manifold $\MCal$ is represented by 
\begin{enumerate}
\item A finite set of players $\NCal = \{1, \ldots, N\}$, with a Riemannian manifold $\MCal_i$ for each player $i \in \NCal$. We denote the action set of player $i$ as $\MCal_i$ and denote the joint action profile as $\MCal = \Pi_{i \in \NCal} \MCal_i$.
\item A loss function $L_i$ is defined for each player $i \in \NCal$, and assigns a real scalar to any joint action profile $x = (x_1, \ldots, x_N) \in \MCal$.
\item An operator $F(\cdot) = (\grad_{x_1} L_1(\cdot), \ldots, \grad_{x_N} L_N(\cdot))$ is defined as the concatenation of gradients of each loss function with respect to its own action.   
\end{enumerate}
\end{definition}
A Riemannian game is smooth if all loss functions are smooth with respect to the metric. The equilibrium notion in a Riemannian game generalizes that in the Euclidean setting: a joint action profile $x^\star \in \MCal$ is a Nash equilibrium if $L_i(x^\star) \leq L_i(x_i, x_{-i}^\star)$ for all players $i \in \NCal$, and all actions $x_i \in \MCal_i$. In the following lemma, we provide a necessary condition for a joint action profile to be a Nash equilibrium of a smooth Riemannian game.
\begin{lemma}\label{Lemma:equilibrium}
The joint action profile $x^\star \in \MCal$ is a Nash equilibrium of a smooth Riemannian game if $\|F(x^\star)\|_{x^\star}= 0$.
\end{lemma}
\begin{definition}\label{def:monotone}
A geodesically $\mu$-strongly monotone game satisfies that for any point $x \in \MCal$, $v \in T_x \MCal$ and one of corresponding geodesics $\gamma$ satisfying $\gamma(0)=x$, $\gamma(1)=y$ and $\dot{\gamma}(0)=v$, we have $\langle \Gamma_y^x F(y)-F(x), v \rangle_x \geq \mu \|v\|_x^2$. If $\mu=0$, we call this game geodesically monotone. 
\end{definition}
Equivalently, a Riemannian game is both geodesically $\mu$-strongly monotone and geodesically $\ell$-smooth if the following statements hold true,  
\begin{equation*}
\inf_{\gamma: [0,1] \mapsto \MCal} \tfrac{\langle \Gamma_{\gamma(1)}^{\gamma(0)}F(\gamma(1)) - F(\gamma(0)), \dot{\gamma}(0)\rangle_{\gamma(0)}}{\|\dot{\gamma}(0)\|_{\gamma(0)}^2} \geq \mu, \quad
\sup_{\gamma: [0,1] \mapsto \MCal} \tfrac{\|F(\gamma(0)) - \Gamma_{\gamma(1)}^{\gamma(0)}F(\gamma(1))\|_{\gamma(0)}}{\|\dot{\gamma}\|_{\gamma(0)}} \leq \ell. 
\end{equation*}
By using the Cauchy-Schwartz inequality, we have
\begin{equation*}
\mu \leq \tfrac{\langle \Gamma_{\gamma(1)}^{\gamma(0)}F(\gamma(1)) - F(\gamma(0)), \dot{\gamma}(0)\rangle_{\gamma(0)}}{\|\dot{\gamma}(0)\|_{\gamma(0)}^2} \leq \tfrac{\|F(\gamma(0)) - \Gamma_{\gamma(1)}^{\gamma(0)}F(\gamma(1))\|_{\gamma(0)}}{\|\dot{\gamma}(0)\|_{\gamma(0)}} \leq \ell. 
\end{equation*}
We now establish two results: (i) the geodesically $\mu$-strongly monotone games are sufficiently general that they cover geodesically $\mu$-strongly convex $\mu$-strongly concave min-max problems; and (ii) the geodesically $\mu$-strongly monotone games are sufficiently structured that $x^\star \in \MCal$ is a Nash equilibrium \textit{if and only if} $\|F(x^\star)\|_{x^\star}= 0$. 
\begin{lemma}\label{Lemma:minimax}
Consider a two-player Riemannian game, where $x_1 \in \MCal_1$ and $x_2 \in \MCal_2$ are the actions of players $1$ and $2$, respectively. The loss of player $1$ is $f(x_1,x_2): \MCal_1 \times \MCal_2 \mapsto \br$ and the loss of player $2$ is $-f(x_1,x_2)$. If $f(x_1, x_2)$ is geodesically $\mu$-strongly convex in $x_1$ and geodesically $\mu$-strongly concave in $x_2$, this game is geodesically $\mu$-strongly monotone.
\end{lemma}
\begin{theorem}\label{Thm:monotone}
The joint action profile $x^\star \in \MCal$ is a Nash equilibrium of a smooth monotone Riemannian game if and only if $\|F(x^\star)\|_{x^\star}= 0$.
\end{theorem}
Recall that $x^\star \in \MCal$ is a Nash equilibrium if $L_i(x^\star) \leq L_i(x_i, x_{-i}^\star)$ for all players $i \in \NCal$, and all actions $x_i \in \MCal_i$. This definition motivates the total gap function in the following form of 
\begin{equation*}
\textnormal{Tgap}(x; D) = \sum_{i \in \NCal} \left(L_i(x) - \min_{z_i \in \BB_{\MCal_i}(x_i; D)} L_i(z_i, x_{-i}) \right), 
\end{equation*}
and guarantees that $x^\star \in \MCal$ is a Nash equilibrium if $\textnormal{Tgap}(x^\star; D)=0$. Another measure of proximity to a Nash equilibrium in Riemannian games is a generalization of the gap function and can be defined as follows, 
\begin{equation*}
\textnormal{gap}(x; D) = \max_{v \in T_x \MCal: \|v\|_x \leq D} \langle F(x), v\rangle_x.  
\end{equation*}
Finally, we establish the relationship between the total gap function, the gap function and the norm of $F$. This explains why we measure the proximity of a joint action profile $x$ to a Nash equilibrium using the norm of $F$ at $x$. 
\begin{lemma}\label{Lemma:relationship}
Suppose that the Riemannian game is geodesically $\mu$-strongly monotone for some $\mu \geq 0$, $x \in \MCal$ and $D>0$. Then, we have
\begin{equation*}
\textnormal{Tgap}(x; D) \leq \textnormal{gap}(x; \sqrt{N}D) \leq \sqrt{N} D \|F(x)\|_x. 
\end{equation*}
\end{lemma}
\begin{remark}
Theorem~\ref{Thm:monotone} is a consequence of Lemmas~\ref{Lemma:equilibrium} and~\ref{Lemma:relationship}. Indeed, by Lemma~\ref{Lemma:equilibrium}, $\|F(x^\star)\|_{x^\star} = 0$ if $x^\star \in \MCal$ is a Nash equilibrium. By Lemma~\ref{Lemma:relationship}, $\textnormal{Tgap}(x; D)=0$ for any $D>0$ if $\|F(x^\star)\|_{x^\star} = 0$, guaranteeing that $x^\star \in \MCal$ is a Nash equilibrium.
\end{remark}

\paragraph{Applications} We provide two examples of Riemannian games to give a better sense of their applicability. One example is a generic model from the optimization and game theory literature~\citep{Monderer-1996-Potential, Sandholm-2001-Potential} and the other is a formalization of a class of applied problems in economics and statistical machine learning~\citep{Pennec-2006-Riemannian, Fletcher-2007-Riemannian}. 
\begin{example}
A game is the Riemannian potential game if there exists a potential function $f: \MCal \mapsto \br$ such that 
\begin{equation*}
L_i(x_i, x_{-i}) - L_i(z_i, x_{-i}) = f(x_i, x_{-i}) - f(z_i, x_{-i}), 
\end{equation*}
for all $i \in \mathcal{N}$, all $x \in \MCal$ and all $z_i \in \MCal_i$. If the function $f$ is geodesically $\mu$-strongly convex, the corresponding Riemannian potential game is geodesically $\mu$-strongly monotone. 
\end{example}
\begin{example}
A robust matrix Karcher mean problem is defined as follows: the Karcher mean of $N$ symmetric positive definite matrices $\{A_i\}_{1 \leq i \leq N}$ is defined as the matrix $X \in \mathcal{M} = \{X \in \br^{n \times n}: X \succ 0, X = X^\top\}$ that minimizes the sum of squared distance induced by the Riemannian metric: $d(X, Y) = \|\textnormal{log}(X^{-1/2}YX^{-1/2})\|_F$. 

We define $f(X; A_1, A_2, \ldots, A_N) = \sum_{i=1}^N (d(X, A_i))^2$ which is nonconvex in the Euclidean setting but is geodesically strongly convex. Then, the robust matrix Karcher mean problem is defined in the form of 
\begin{equation*}
\min_{X \in \MCal} \max_{Y_i \in \MCal} \ f(X; Y_1, Y_2, \ldots, Y_N) - \gamma\left(\sum_{i=1}^N (d(Y_i, A_i))^2\right), 
\end{equation*}
where $\gamma > 0$ stands for the trade-off between the computation of Karcher mean over $\{Y_i\}_{i=1}^N$ and the difference between $\{A_i\}_{i=1}^N$ and $\{Y_i\}_{i=1}^N$. This problem is geodesically strongly convex in $X$ and geodesically strongly concave in $(Y_1, Y_2, \ldots, Y_N)$ when $\gamma$ is sufficiently large. 
\end{example}
\begin{example}
In network communication, nodes often make decisions influenced by the underlying geometry of a data manifold. When each node's outcome depends on the actions of its neighbors, such interactions yield a game on manifolds. Riemannian games thus provide a natural framework for predicting and analyzing equilibrium behaviors, particularly in networks whose topology or data flow exhibits geometric complexity.
\end{example}
\begin{example}
In robotics and control, especially within multi-robot systems, the decision or configuration spaces are frequently modeled as manifolds. Riemannian games yield valuable insights into optimal navigation, coordination, and control strategies when agents' actions are interdependent, enabling a principled understanding of cooperative or competitive behaviors.
\end{example}
\begin{example}
In quantum systems, the set of quantum states forms a manifold. When decision-making or computation among quantum entities involves strategic interdependence, the framework of Riemannian games offers a powerful perspective for studying equilibrium phenomena. This approach can aid in the design of quantum algorithms and the analysis of multi-particle interactions governed by geometric constraints.
\end{example}
Riemannian games subsume all general games defined in Euclidean setting as well as both minimization and maximization problems posed on geodesic spaces. A representative class of such problems is min-max optimization on geodesic manifolds, which serves as an abstraction for many machine learning applications, including principal component analysis~\citep{Boumal-2011-RTRMC}, dictionary learning~\citep{Sun-2016-Complete-I, Sun-2016-Complete-II}, deep neural networks~\citep{Huang-2018-Orthogonal}, and low-rank matrix learning~\citep{Vandereycken-2013-Low, Jawanpuria-2018-Unified}.
Indeed, principal component analysis can be formulated as an optimization problem on the Grassmann manifold, which captures the subspace geometry inherent in the data.

%% file: sec/algorithm.tex
\section{Geometry-Agnostic Convergence Rates} \label{sec:results}
We prove that deterministic and stochastic Riemannian gradient descent methods attain a linear, geometry-agnostic, and last-iterate convergence rate in geodesically smooth and strongly monotone games. We further extend these results to fully adaptive Riemannian gradient descent ascent methods that require no prior knowledge of problem parameters.

\subsection{Deterministic and Stochastic Riemannian Gradient Descent}
The assumptions summarized in Assumption~\ref{Assumption:main-1} guarantee the existence of at least one Nash equilibrium in this class of Riemannian games and the linear, geometry-agnostic, and last-iterate convergence rate of Algorithm~\ref{alg:SRGD} for finding one Nash equilibrium. 
\begin{assumption}\label{Assumption:main-1}
We make the following assumptions:
\begin{enumerate}
\item The game is geodescially $\mu$-strongly monotone and $\ell$-smooth with known $\mu$ and $\ell$. 
\item The manifold $\MCal$ does not have boundaries and is complete under the metric $d_\MCal(\cdot, \cdot)$. 
\item In the stochastic setting, we have $d_{\MCal}(x_0, x^\star) \leq B$ for the initial joint action $x_0$.
\end{enumerate}
\end{assumption}
\begin{algorithm}[!t] 
\begin{algorithmic}
\caption{Deterministic and Stochastic Riemannian Gradient Descent} \label{alg:SRGD}
\STATE \textbf{Input:} initial profile $x^0$, stepsizes $\eta_t$, minibatch sizes $m_t$, and iteration number $T \geq 1$. 
\FOR{$t = 0,1,\ldots, T-1$}
\STATE Query $g_i^t$ as a noisy estimator of $F(x^t) = (\grad_{x_1} L_1(x^t), \ldots, \grad_{x_N} L_N(x^t))$ with mean $0$ and variance $\sigma^2>0$ for $1 \leq i \leq m_t$. 
\STATE $x^{t+1} \gets \Exp_{x^t}(-\eta_t \cdot g^t)$ where $g^t = \frac{1}{m_t}(\sum_{i=1}^{m_t} g_i^t)$. 
\ENDFOR
\STATE \textbf{Return:} $x^T$.
\end{algorithmic}
\end{algorithm}
We establish a descent inequality to analyze Algorithm~\ref{alg:SRGD} and its fully adaptive variant, where $\hat{\mu}$ and $\hat{\ell}$ are sufficiently accurate estimates of strong monotonicity and smoothness parameters. For Algorithm~\ref{alg:SRGD}, these estimates can be regarded as the exact parameter values under which our condition holds, without loss of generality or intuition.
\begin{lemma}\label{Lemma:descent}
Under Assumption~\ref{Assumption:main-1} and let $\{x^t\}_{t=0}^T$ be generated by Algorithm~\ref{alg:SRGD} with stepsizes $\eta_t > 0$ and minibatch sizes $m_t \geq 1$. Then, we have 
\begin{equation*}
\EE[\|F(x^{t+1})\|_{x^{t+1}}^2] \leq \left(1-\tfrac{2\eta_t\hat{\mu}-(\eta_t\hat{\ell})^2}{2-2\eta_t\hat{\mu}+(\eta_t\hat{\ell})^2}\right)\EE[\|F(x^t)\|^2_{x^t}] + \tfrac{4\sigma^2}{m_t\eta_t(2\hat{\mu}-\eta_t\hat{\ell}^2)}, 
\end{equation*}
where $\hat{\mu}$, $\hat{\ell}$ and $\eta_t$ satisfy 
\begin{equation*}
0 < \eta_t < \tfrac{2\hat{\mu}}{\hat{\ell}^2}, \quad 0 < \hat{\mu} \leq \tfrac{\EE[\langle \Gamma_{x^{t+1}}^{x^t} F(x^{t+1}) - F(x^t), - \eta_t g^t\rangle_{x^t}]}{\EE[\|\eta_t g^t\|_{x^t}^2]} \leq \sqrt{\tfrac{\EE[\|\Gamma_{x^{t+1}}^{x^t} F(x^{t+1}) - F(x^t)\|_{x^t}^2]}{\EE[\|\eta_t g^t\|_{x^t}^2]}} \leq \hat{\ell}. 
\end{equation*}
\end{lemma}
\begin{proof}
Since $g_i^t$ is a noisy estimator of $F(x^t) = (\grad_{x_1} L_1(x^t), \ldots, \grad_{x_N} L_N(x^t))$ with mean $0$ and variance $\sigma^2>0$ for $1 \leq i \leq m_t$, and $g^t = \frac{1}{m_t}(\sum_{i=1}^{m_t} g_i^t)$, we have 
\begin{equation*}
\EE[g^t \mid x^t] = F(x^t), \quad \EE[\|g^t - F(x^t)\|^2_{x^t} \mid x^t] \leq \tfrac{\sigma^2}{m_t}. 
\end{equation*}
This implies 
\begin{equation*}
\EE[\|F(x^t)\|^2_{x^t}] - \EE[\langle \Gamma_{x^{t+1}}^{x^k} F(x^{t+1}), g^t \rangle_{x^t}] = \EE[\langle \Gamma_{x^{t+1}}^{x^t} F(x^{t+1})-F(x^t),  -g^t \rangle_{x^t}].  
\end{equation*}
By the definition of $\hat{\mu}$, we have 
\begin{equation*}
\EE[\langle \Gamma_{x^{t+1}}^{x^t} F(x^{t+1})-F(x^t),  -\eta_t \cdot g^t \rangle_{x^t}] \geq \hat{\mu}(\eta_t)^2 \left(\EE[\|g^t\|_{x^t}^2]\right). 
\end{equation*}
Putting these pieces together yields 
\begin{equation}\label{inequality:descent-1}
\EE[\|F(x^t)\|^2_{x^t}] - \EE[\langle \Gamma_{x^{t+1}}^{x^k} F(x^{t+1}), g^t \rangle_{x^t}] \geq \tfrac{1}{\eta_t} \cdot \hat{\mu}(\eta_t)^2\left(\EE[\|g^t\|_{x^t}^2]\right) = \hat{\mu}\eta_t\left(\EE[\|g^t\|_{x^t}^2]\right). 
\end{equation}
Since $\EE[g^t \mid x^t] = F(x^t)$, we have 
\begin{equation*}
\EE[\|g^t\|_{x^t}^2 \mid x^t] = \|F(x^t)\|_{x^t}^2 + \EE[\|g^t - F(x^t)\|^2_{x^t} \mid x^t]. 
\end{equation*}
By taking the expectation of both sides of the above equation and combining it with Eq.~\eqref{inequality:descent-1}, we have 
\begin{equation*}
\EE[\|F(x^t)\|^2_{x^t}] - \EE[\langle \Gamma_{x^{t+1}}^{x^k} F(x^{t+1}), g^t \rangle_{x^t}] \geq \hat{\mu} \eta_t \left(\EE[\|F(x^t)\|_{x^t}^2] + \EE[\|g^t - F(x^t)\|^2_{x^t}]\right). 
\end{equation*}
Equivalently, we have
\begin{equation}\label{inequality:descent-2}
(1 - \hat{\mu}\eta_t)\EE[\|F(x^t)\|^2_{x^t}] - \hat{\mu}\eta_t\EE[\|g^t-F(x^t)\|^2_{x^t}] - \EE[\langle \Gamma_{x^{t+1}}^{x^t} F(x^{t+1}), g^t\rangle_{x^t}] \geq 0. 
\end{equation}
By the definition of $\hat{\ell}$, we have 
\begin{equation*}
\EE[\|\Gamma_{x^{t+1}}^{x^t} F(x^{t+1}) - F(x^t)\|_{x^t}^2] \leq (\hat{\ell})^2 \EE[\|\eta_t g^t\|^2_{x^t}] = (\ell\eta_t)^2 \EE[\|g^t\|^2_{x^t}]. 
\end{equation*}
Note that $\EE[\|g^t\|^2_{x^t}] = \EE[\|F(x^t)\|_{x^t}^2] + \EE[\|g^t - F(x^t)\|^2_{x^t}]$. Then, we have 
\begin{equation*}
\EE[\|\Gamma_{x^{t+1}}^{x^t} F(x^{t+1}) - F(x^t)\|_{x^t}^2] \leq (\ell\eta_t)^2 \left(\EE[\|F(x^t)\|_{x^t}^2] + \EE[\|g^t - F(x^t)\|^2_{x^t}]\right). 
\end{equation*}
By definition, we have 
\begin{eqnarray*}
\lefteqn{\|F(x^{t+1})\|_{x^{t+1}}^2 = \langle F(x^{t+1}),F(x^{t+1}) \rangle_{x^{t+1}}} \\ 
& = & \langle \Gamma_{x^{t+1}}^{x^t}F(x^{t+1}),\Gamma_{x_{t+1}}^{x^t} F(x^{t+1}) \rangle_{x^t} = \|\Gamma_{x^{t+1}}^{x^t} F(x^{t+1})\|^2_{x^t} \nonumber. 
\end{eqnarray*}
This implies 
\begin{eqnarray*}
\EE[\|\Gamma_{x^{t+1}}^{x^t} F(x^{t+1}) - F(x^t)\|_{x^t}^2] & = & \EE[\|F(x^t)\|^2_{x^t}] - 2\EE[\langle \Gamma_{x^{t+1}}^{x^t} F(x^{t+1}), F(x^t) \rangle_{x^t}] + \EE[ \|\Gamma_{x^{t+1}}^{x^t} F(x^{t+1})\|^2_{x^t}] \\
& = & \EE[\|F(x^t)\|^2_{x^t}] - 2\EE[\langle \Gamma_{x^{t+1}}^{x^t} F(x^{t+1}), F(x^t) \rangle_{x^t}] + \EE[\|F(x^{t+1})\|_{x^{t+1}}^2]. 
\end{eqnarray*}
Putting these pieces together yields 
\begin{eqnarray}\label{inequality:descent-3}
\MoveEqLeft{{((\eta_t\hat{\ell})^2-1)\EE[\|F(x^t)\|_{x^t}^2] + (\eta_t\hat{\ell})^2 \EE[\|g^t-F(x^t)\|_{x^t}^2]}} \\
 & + 2\EE[\langle \Gamma_{x^{t+1}}^{x^t} F(x^{t+1}), F(x^t) \rangle_{x^t}] \geq \EE[\|F(x^{t+1})\|^2_{x^{t+1}}]. \nonumber
\end{eqnarray}
Combining Eq.~\eqref{inequality:descent-2} and Eq.~\eqref{inequality:descent-3} yields 
\begin{eqnarray} \label{inequality:descent-4}
\MoveEqLeft{(1-2\eta_t\hat{\mu} + (\eta_t\hat{\ell})^2)\EE[\|F(x^t)\|^2_{x^t}] + ((\eta_t\hat{\ell})^2-2\hat{\mu}\eta_t)\EE[\|g^t - F(x^t)\|^2_{x^t}]}  \\
& + 2\EE[\langle \Gamma_{x^{t+1}}^{x^t} F(x^{t+1}), F(x^t)-g^t\rangle_{x^t}] \geq \EE[\|F(x^{t+1})\|^2_{x^{t+1}}]. \nonumber
\end{eqnarray}
Since $0 < \eta_t < \frac{2\hat{\mu}}{\hat{\ell}^2}$ and $0 < \hat{\mu} \leq \hat{\ell}$, we have 
\begin{equation*}
(\eta_t\hat{\ell})^2 -2\hat{\mu}\eta_t \leq 0, \quad 2 - 2\hat{\mu}\eta_t + (\eta_t\hat{\ell})^2 \geq 2 - 2\eta_t\hat{\mu} + (\eta_t\hat{\mu})^2 \geq 1. 
\end{equation*}
By using Young's inequality and the property of parallel transport, we have 
\begin{eqnarray*}
\EE[\langle \Gamma_{x^{t+1}}^{x^t} F(x^{t+1}), F(x^t)-g^t\rangle_{x^t}] & \leq & \tfrac{2\hat{\mu}\eta_t-(\eta_t\hat{\ell})^2}{4}\EE[\|g^t - F(x^t)\|^2_{x^t}] + \tfrac{1}{2\hat{\mu}\eta_t-(\eta_t\hat{\ell})^2}\EE[\|\Gamma_{x^{t+1}}^{x^t} F(x^{t+1})\|^2_{x^t}] \\ 
& = & \tfrac{2\hat{\mu}\eta_t-(\eta_t\hat{\ell})^2}{4}\EE[\|g^t - F(x^t)\|^2_{x^t}] + \tfrac{1}{2\hat{\mu}\eta_t-(\eta_t\hat{\ell})^2}\EE[\|F(x^{t+1})\|^2_{x^{t+1}}]. 
\end{eqnarray*}
In addition, we have 
\begin{equation*}
\EE[\|g^t - F(x^t)\|^2_{x^t}] = \EE\left[\EE[\|g^t - F(x^t)\|^2_{x^t} \mid x^t]\right] \leq \tfrac{\sigma^2}{m_t}. 
\end{equation*}
Putting these pieces together with Eq.~\eqref{inequality:descent-4} yields 
\begin{equation*}
(1-2\hat{\mu}\eta_t + (\eta_t\hat{\ell})^2)\EE[\|F(x^t)\|^2_{x^t}] + \tfrac{2\sigma^2}{m_t\eta_t(2\hat{\mu}-\eta_t \hat{\ell}^2)} \geq \left(1-\tfrac{2\hat{\mu}\eta_t-(\eta_t\hat{\ell})^2}{2}\right)\EE[\|F(x^{t+1})\|^2_{x^{t+1}}]. 
\end{equation*}
Equivalently, we have 
\begin{equation*}
\left(1-\tfrac{2\hat{\mu}\eta_t-(\eta_t\hat{\ell})^2}{2-2\hat{\mu}\eta_t+(\eta_t\hat{\ell})^2}\right)\EE[\| F(x^t)\|^2_{x^t}] + \tfrac{4\sigma^2}{m_t\eta_t(2\hat{\mu}-\eta_t\hat{\ell}^2)(2-2\hat{\mu}\eta_t+(\eta_t \hat{\ell})^2)} \geq \EE[\|F(x^{t+1})\|^2_{x^{t+1}}]. 
\end{equation*}
Since $2 - 2\hat{\mu}\eta_t + (\eta_t\hat{\ell})^2 \geq 1$, we have 
\begin{equation*}
\left(1-\tfrac{2\hat{\mu}\eta_t-(\eta_t\hat{\ell})^2}{2-2\hat{\mu}\eta_t+(\eta_t\hat{\ell})^2}\right)\EE[\| F(x^t)\|^2_{x^t}] + \tfrac{4\sigma^2}{m_t\eta_t(2\hat{\mu}-\eta_t\hat{\ell}^2)} \geq \EE[\|F(x^{t+1})\|^2_{x^{t+1}}]. 
\end{equation*}
This completes the proof. 
\end{proof}
\begin{corollary}\label{Corollary:exact}
If the strong monotonicity and smoothness parameters $\mu$ and $\ell$ are known, we choose $\hat{\mu}=\mu$ and $\hat{\ell}=\ell$ such that the conditions in Lemma~\ref{Lemma:descent} are satisfied for all $t \geq 1$.  
\end{corollary}
\begin{proof}
We consider the geodesic with $\gamma(0)=x^t$, $\gamma(1)=x^{t+1}$ and $\dot{\gamma}(0)=-\eta_t g^t$. Since the game is geodescially $\mu$-strongly monotone and $\ell$-smooth, we have 
\begin{eqnarray*}
\langle \Gamma_{x^{t+1}}^{x^t}F(x^{t+1})-F(x^t), -\eta_t g^t \rangle_{x^t} & \geq & \mu\|\eta_t g^t\|_{x^t}^2, \\ 
\|F(x^t) - \Gamma_{x^{t+1}}^{x^t} F(x^{t+1})\|_{x^t} & \leq & \ell\|\eta_t g^t\|_{x^t}. 
\end{eqnarray*}
This implies that $\hat{\mu}=\mu$ and $\hat{\ell}=\ell$ satisfy the conditions in Lemma~\ref{Lemma:descent} for all $t \geq 1$.  
\end{proof}
We prove the existence of at least one Nash equilibrium using Lemmas~\ref{Lemma:relationship} and~\ref{Lemma:descent}, and present the convergence guarantees for Algorithm~\ref{alg:SRGD} in Theorem~\ref{Thm:SRGD}. 
\begin{lemma}\label{Lemma:existence}
A Riemannian game has at least one Nash equilibrium under Assumption~\ref{Assumption:main-1}.
\end{lemma}
\begin{proof}
We implement Algorithm~\ref{alg:SRGD} with $g_i^t = F(x^t)$ for all $i$, a constant stepsize $\eta_t = \frac{\mu}{\ell^2}$ and an initial profile $x^0 \in \MCal$. By choosing $\hat{\mu} = \mu$ and $\hat{\ell} = \ell$ and defining $\kappa = \frac{\ell}{\mu}$, Lemma~\ref{Lemma:descent} implies 
\begin{equation*}
\|F(x^t)\|_{x^t}^2 \leq \left(1 - \tfrac{1}{2\kappa^2-1}\right)^t \|F(x^0)\|^2_{x^0} \leq e^{-\frac{t}{2\kappa^2}}\|F(x^0)\|^2_{x^0}. 
\end{equation*}
By definition, we have $d_\MCal(x^{t+1}, x^t) \leq \eta_t\|F(x^t)\|_{x^t}$. Thus, we have 
\begin{equation*}
\sum_{k\geq 1} d_\MCal(x^{t+1}, x^t) \leq \left(\sum_{k\geq 1} e^{-\frac{k}{2\kappa^2}}\right)\tfrac{1}{\kappa\ell}\|F(x^0)\|^2_{x^0}. 
\end{equation*}
Since the right-hand side of the above inequality is bounded by $D_0$ and $(\MCal, d_\MCal)$ is complete, the sequence of iterates $\{x^t\}_{t \geq 1}$ is Cauchy and converges to $x^\star = \lim_{t \rightarrow +\infty} x^t \in \BB_\MCal(x^0; D_0)$. In addition, we have 
\begin{equation*}
\|F(x^\star)\|_{x^\star} = \lim_{t \rightarrow +\infty} \|F(x^t)\|_{x^t} \leq \left(\lim_{t \rightarrow +\infty} e^{-\frac{t}{2\kappa^2}}\right) \|F(x^0)\|^2_{x^0} = 0. 
\end{equation*}
This together with Lemma~\ref{Lemma:equilibrium} yields the desired result. 
\end{proof}

\begin{theorem}\label{Thm:SRGD}
Under Assumption~\ref{Assumption:main-1} and let $\{x^t\}_{t=0}^T$ be generated by Algorithm~\ref{alg:SRGD} with stepsizes $\eta_t = \frac{\mu}{\ell^2}$ for all $t \geq 0$. Then, we have 
\begin{itemize}
\item In the deterministic setting with $m_t = 1$ for all $t$, there exists some $T>0$ such that the output of Algorithm~\ref{alg:SRGD} satisfies $\|F(x^T)\|_{x^T} \leq \epsilon$ and the number of Riemannian gradient evaluations is bounded by
\begin{equation*}
O\left(\kappa^2\log\left(\frac{\ell d_\MCal(x^0,x^\star)}{\epsilon}\right)\right). 
\end{equation*}
\item In the stochastic setting with $m_t = \lceil \frac{16\kappa^4\sigma^2}{\ell^2B^2}e^{\frac{t}{4\kappa^2}} \rceil$ for all $t$, there exists some $T>0$ such that the output of Algorithm~\ref{alg:SRGD} satisfies $\EE[\|F(x^T)\|_{x^T}] \leq \epsilon$ and the number of noisy Riemannian gradient evaluations is bounded by
\begin{equation*}
O\left(\kappa^2\log\left(\frac{\ell B}{\epsilon}\right) + \frac{\kappa^6\sigma^2}{\epsilon^2 }\right). 
\end{equation*}
\end{itemize}
In the above results, $\kappa = \ell/\mu \geq 1$ is defined as the condition number of the game. 
\end{theorem}
\begin{proof}
By Lemma~\ref{Lemma:equilibrium}, we have $F(x^\star)=0$. Since the game is geodescially $\ell$-smooth, we consider the distance-minimizing constant-speed geodesic with $\gamma(0)=x^0$, $\gamma(1)=x^\star$ and $\dot{\gamma}(0)=v$ and have 
\begin{equation*}
\|F(x^0)\|_{x^0} = \|\Gamma_{x^0}^{x^\star} F(x^0)\|_{x^\star} = \|\Gamma_{x^0}^{x^\star} F(x^0) - F(x^\star)\|_{x^\star} \leq \ell\|v\|_{x^\star} = \ell d_\MCal(x^0,x^\star).
\end{equation*}
Under Assumption~\ref{Assumption:main-1}, we obtain from Corollary~\ref{Corollary:exact} that $\hat{\mu}=\mu$ and $\hat{\ell}=\ell$. Since $\eta_t = \frac{\mu}{\ell^2}$ for all $t \geq 0$, Lemma~\ref{Lemma:descent} implies
\begin{equation*}
\EE[\|F(x^{t+1})\|_{x^{t+1}}^2] \leq \left(1-\tfrac{1}{2\kappa^2-1}\right)\EE[\|F(x^t)\|^2_{x^t}] + \tfrac{4\kappa^2\sigma^2}{m_t} \leq \left(1-\tfrac{1}{2\kappa^2}\right)\EE[\|F(x^t)\|^2_{x^t}] + \tfrac{4\kappa^2\sigma^2}{m_t}. 
\end{equation*}
In the deterministic setting with $m_t=1$, we have $\sigma^2=0$ and no randomness. Thus, we have
\begin{equation*}
\|F(x^T)\|_{x^T}^2 \leq \left(1-\tfrac{1}{2\kappa^2}\right)^T \|F(x^0)\|_{x^0}^2 \leq e^{-\frac{T}{2\kappa^2}} \ell^2(d_\MCal(x^0,x^\star))^2. 
\end{equation*}
This implies that the number of Riemannian gradient evaluations required by Algorithm~\ref{alg:SRGD} to return a point $x^T$ satisfying $\|F(x^T)\|_{x^T} \leq \epsilon$ is bounded by 
\begin{equation*}
O\left(\kappa^2\log\left(\frac{\ell d_\MCal(x^0,x^\star)}{\epsilon}\right)\right). 
\end{equation*}
In the stochastic setting with $m_t = \lceil \frac{16\kappa^4\sigma^2}{\ell^2B^2}e^{\frac{t}{4\kappa^2}} \rceil$, we have 
\begin{equation*}
\EE[\|F(x^{t+1})\|_{x^{t+1}}^2] \leq \left(1-\tfrac{1}{2\kappa^2}\right)\EE[\|F(x^t)\|^2_{x^t}] + \tfrac{\ell^2 B^2}{4\kappa^2}e^{-\frac{t}{4\kappa^2}}. 
\end{equation*}
In what follows, we prove that $\EE[\|F(x^t)\|^2_{x^t}] \leq \ell^2 B^2 e^{-\frac{t}{4\kappa^2}}$. Since $\|F(x^0)\|^2_{x^0} \leq \ell B$, we have that this inequality holds true if $t=0$. By using inductive arguments, we have 
\begin{equation*}
\EE[\|F(x^{t+1})\|_{x^{t+1}}^2] \leq \left(1-\tfrac{1}{2\kappa^2}\right)\ell^2 B^2 e^{-\frac{t}{4\kappa^2}} + \tfrac{\ell^2 B^2}{4\kappa^2}e^{-\frac{t}{4\kappa^2}} = \left(1-\tfrac{1}{4\kappa^2}\right)\ell^2 B^2 e^{-\frac{t}{4\kappa^2}} \leq \ell^2 B^2 e^{-\frac{t+1}{4\kappa^2}}. 
\end{equation*}
This implies that the number of iterations required by Algorithm~\ref{alg:SRGD} to return a point $x^T$ satisfying $\EE[\|F(x^T)\|_{x^T}] \leq \epsilon$ is bounded by 
\begin{equation*}
O\left(\kappa^2\log\left(\frac{\ell B}{\epsilon}\right)\right). 
\end{equation*}
As a consequence, the number of noisy Riemannian gradient evaluations required by Algorithm~\ref{alg:SRGD} to return a point $x^T$ satisfying $\EE[\|F(x^T)\|_{x^T}] \leq \epsilon$ is bounded by 
\begin{equation*}
\sum_{t=0}^{T-1} m_t \leq T + \tfrac{16\kappa^4\sigma^2}{\ell^2 B^2} \left(\sum_{t=0}^{T-1} e^{\frac{t}{4\kappa^2}}\right) = T + \tfrac{16\kappa^4\sigma^2}{\ell^2 B^2} \left(\tfrac{e^{\frac{T}{4\kappa^2}}-1}{e^{\frac{1}{4\kappa^2}}-1}\right). 
\end{equation*}
Since $T=O(\kappa^2\log(\frac{\ell B}{\epsilon}))$ and $e^{\frac{1}{4\kappa^2}}-1 \geq \frac{1}{4\kappa^2}$, we have 
\begin{equation*}
\sum_{t=0}^{T-1} m_t = O\left(\kappa^2\log\left(\frac{\ell B}{\epsilon}\right) + \frac{\kappa^6\sigma^2}{\epsilon^2 }\right). 
\end{equation*}
This completes the proof. 
\end{proof}
\begin{remark}
Theorem~\ref{Thm:SRGD} illustrates the geometry-agnostic and last-iterate convergence of Algorithm~\ref{alg:SRGD} for solving geodesically strongly monotone Riemannian games. Within the context of Riemannian min-max optimization, our results can be directly compared with those obtained for the Riemannian corrected extragradient method (RCEG) in~\citet{Jordan-2022-First}. In particular, their convergence bound is expressed in terms of $d_\MCal(x^T, x^\star)$, which relates to our measure through the inequality $\|F(x^T)\|_{x^T} \leq \ell d_\MCal(x^T, x^\star)$. The key advantage of our analysis is that both the stepsize selection and the convergence rate are entirely independent of any \textit{sectional curvature} parameters, thereby demonstrating the robustness and generality of our method.
\end{remark}
\begin{remark}
Our analysis extends to locally geodesically strongly monotone settings. The key lies in how the local region is defined. Indeed, if a game is geodesically strongly monotone within $\{x \in \MCal: d_{\MCal_i}(x_i, x_i^\star) \leq \delta \textnormal{ for all } i \in \NCal\}$, then any smaller neighborhood $\{x \in \MCal: d_\MCal(x, x^\star) \leq \delta\}$ is contained within such a region, and the game remains geodesically strongly monotone therein. When the initialization $x^0$ lies in this neighborhood, our theoretical analysis guarantees a last-iterate linear convergence rate. This argument and the associated notion of a local region are standard in the analysis of min-max optimization problems in both Euclidean and Riemannian settings~\citep{Liang-2019-Interaction, Jordan-2022-First}. As an illustrative example of a globally geodesically strongly monotone game, we highlight the robust matrix Karcher mean problem.
\end{remark}

\subsection{Fully Adaptive Riemannian Gradient Descent}\label{sec:FARGD}
The assumptions summarized in Assumption~\ref{Assumption:main-2} guarantee the linear, geometry-agnostic, and last-iterate convergence rate of Algorithm~\ref{alg:FARGD} for finding one Nash equilibrium. 
\begin{algorithm}[!t]
\caption{Fully Adaptive Riemannian Gradient Descent}
\label{alg:FARGD}
\begin{algorithmic}
\STATE \textbf{Input:} initial joint action $x^0$, initial estimates of strong monotonicity and smoothness $\hat{\mu}_0$ and $\hat{\ell}_0$, and iteration number $T \geq 1$. 
\FOR{$t = 0,1,\ldots, T-1$}
\STATE Query an exact quantity of $F(x^t) = (\grad_{x_1} L_1(x^t), \ldots, \grad_{x_N} L_N(x^t))$. 
\STATE $\tilde{x}^{t+1} \gets \Exp_{x^t}(-\eta_t \cdot F(x^t))$ where $\eta_t = \frac{\hat{\mu}_t}{\hat{\ell}_t^2}$. 
\IF{$\hat{\mu}_t > \frac{\langle \Gamma_{\tilde{x}^{t+1}}^{x^t} F(\tilde{x}^{t+1}) - F(x^t), -\eta_t \cdot F(x^t) \rangle_{x^t}}{\|\eta_t F(x^t)\|^2_{x^t}}$}
\STATE $\hat{\mu}_t \gets \frac{\hat{\mu}_t}{2}$ and restart iteration $t$. 
\ELSIF{$\hat{\ell}_t < \frac{\| \Gamma_{\tilde{x}^{t+1}}^{x^t} F(\tilde{x}^{t+1}) - F(x^t)\|_{x^t}}{\|\eta_t F(x^t)\|_{x^t}}$}
\STATE $\hat{\ell}_t \gets 2\hat{\ell}_t$ and restart iteration $t$. 
\ELSE
\STATE $\hat{\mu}_{t+1} \gets \hat{\mu}_t$, $\hat{\ell}_{t+1} \gets \hat{\ell}_t$ and $x^{t+1} \gets \tilde{x}^{t+1}$. 
\ENDIF
\ENDFOR
\STATE \textbf{Return:} $x^T$.
\end{algorithmic}
\end{algorithm}
\begin{assumption}\label{Assumption:main-2}
We make the following assumptions:
\begin{itemize}
\item The game is geodescially $\mu$-strongly monotone and $\ell$-smooth with unknown $\mu$ and $\ell$. 
\item The manifold $\MCal$ does not have boundaries and is complete under the metric $d_\MCal(\cdot, \cdot)$. 
\item The initial estimates for strong monotonicity and smoothness parameters $\hat{\mu}_0$ and $\hat{\ell}_0$ satisfy that $\mu < \hat{\mu}_0 \leq \hat{\ell}_0 < \ell$.
\end{itemize}
\end{assumption}
Identifying values that satisfy the last item in Assumption~\ref{Assumption:main-2} can be achieved as follows: we query a point $x \in \MCal$ and $v \in T_x \MCal$, and construct a constant-speed geodesic $\gamma$ with $\gamma(0) = x$, $\gamma(1) = y$ and $\dot{\gamma}(0) = v$. Then, we have 
\begin{equation*}
\hat{\mu}_0 = \tfrac{\langle \Gamma_y^x F(y) - F(x), v \rangle_x}{\|v\|_x^2}, \quad \hat{\ell}_0= \tfrac{\|\Gamma_y^x F(y)-F(x)\|_x}{\|v\|_x}.
\end{equation*}
Our algorithm incorporates two additional inputs: an initial estimate of strong monotonicity $\hat{\mu}_0$, and an initial estimate of smoothness $\hat{\ell}_0$. The structure remains similar to Algorithm~\ref{alg:SRGD} while the key difference lies in its adaptive scheme. In particular, whenever an iterate fails to satisfy the $\hat{\mu}$-strong monotonicity condition, the algorithm halves the current estimate of $\hat{\mu}$ and repeats the iteration. If an iterate violates the $\hat{\ell}$-smoothness condition, the algorithm doubles the current estimate of $\hat{\ell}$ and repeats the iteration. 
\begin{theorem}\label{Thm:FARGD}
Under Assumption~\ref{Assumption:main-2} and let $\{x^t\}_{t=0}^T$ be generated by Algorithm~\ref{alg:FARGD} with stepsizes $\eta_t = \frac{\hat{\mu}_t}{\hat{\ell}_t^2}$ for all $t \geq 0$. There exists some $T>0$ such that the output of Algorithm~\ref{alg:FARGD} satisfies $\|F(x^T)\|_{x^T} \leq \epsilon$ and the number of Riemannian gradient evaluations is bounded by
\begin{equation*}
O\left(\kappa^2\left(\log\left(\frac{\mu}{\hat{\mu}_0}\right)+\log\left(\frac{\widehat{\ell}_0}{\ell}\right)\right)\log\left(\frac{\ell d_\MCal(x^0,x^\star)}{\epsilon}\right)\right). 
\end{equation*}
In this result, $\kappa = \ell/\mu \geq 1$ is defined as the condition number of the game. 
\end{theorem}
\begin{proof}
We first bound the number of repeated iterations. Indeed, Algorithm~\ref{alg:FARGD} restarts the iteration either when we half the value of $\hat{\mu}$ or double the value of $\hat{\ell}$. In the former case, when $\hat{\mu}_t \leq \mu$, the following condition will never be true, 
\begin{equation*}
\hat{\mu}_t > \frac{\langle \Gamma_{\tilde{x}^{t+1}}^{x^t} F(\tilde{x}^{t+1}) - F(x^t)), -\eta_t \cdot F(x^t) \rangle_{x^t}}{\|\eta_t F(x^t)\|^2_{x^t}}. 
\end{equation*}
Thus, we never half the value of $\hat{\mu}_t$ after $\hat{\mu}_t \leq \mu$ and $\hat{\mu}_t \geq \frac{\mu}{2}$ for all $t \geq 0$. In addition, the number of repeated iterations for updating $\hat{\mu}$ is bounded by $\lceil \log(\frac{\mu}{\hat{\mu}_0})\rceil$. The similar argument holds for the latter case. Indeed, $\hat{\ell}_t \leq 2\ell$ for all $t \geq 0$ and the number of repeated iterations for updating $\hat{\ell}$ is bounded by $\lceil\log(\frac{\hat{\ell}_0}{\ell})\rceil$. 

It suffices to quantify the progress during the iterations where we do not repeat. Indeed, the following two conditions hold true,  
\begin{equation*}
\hat{\mu} \leq \tfrac{\langle \Gamma_{x^{t+1}}^{x^t} F(x^{t+1}) - F(x^t), -\eta_t \cdot F(x^t) \rangle_{x^t}}{\|\eta_t F(x^t)\|^2_{x^t}}, \quad \hat{\ell} \geq \tfrac{\|\Gamma_{x^{t+1}}^{x^t} F(x^{t+1}) - F(x^t)\|_{x^t}}{\|\eta_t F(x^t)\|_{x^t}}. 
\end{equation*}
By Assumption~\ref{Assumption:main-2}, we have $\hat{\mu}_0 \leq \hat{\ell}_0$. This together with $\hat{\mu}_t \leq \hat{\mu}_0$ and $\hat{\ell}_t \geq \hat{\ell}_0$ for all $t \geq 0$ yields $\hat{\mu}_t \leq \hat{\ell}_t$ for all $t \geq 0$. In addition, $\eta_t = \frac{\hat{\mu}_t}{\hat{\ell}_t^2}$ for all $t \geq 0$. Then, Lemma~\ref{Lemma:descent} implies 
\begin{equation*}
\|F(x^{t+1})\|_{x^{t+1}}^2 \leq \left(1-\tfrac{1}{2\left(\frac{\hat{\ell}_t}{\hat{\mu}_t}\right)^2-1}\right)\|F(x^t)\|^2_{x^t}. 
\end{equation*}
Since $\hat{\mu}_t \geq \frac{\mu}{2}$ and $\hat{\ell}_t \leq 2\ell$ for all $t \geq 0$, we have $1 \leq \frac{\hat{\ell}_t}{\hat{\mu}_t} \leq 4\kappa$. Then, we have 
\begin{equation*}
\|F(x^{t+1})\|_{x^{t+1}}^2 \leq \left(1-\tfrac{1}{32\kappa^2-1}\right)\|F(x^t)\|^2_{x^t} \leq \left(1-\tfrac{1}{32\kappa^2}\right)\|F(x^t)\|^2_{x^t}. 
\end{equation*}
Thus, after $T$ iterations where we do not repeat, we have 
\begin{equation*}
\|F(x^T)\|_{x^T}^2 \leq \left(1-\tfrac{1}{32\kappa^2}\right)^T \|F(x^0)\|_{x^0}^2 \leq e^{-\frac{T}{32\kappa^2}} \ell^2(d_\MCal(x^0,x^\star))^2. 
\end{equation*}     
This implies that the number of iterations where we do not repeat required by Algorithm~\ref{alg:FARGD} to return a point $x^T$ satisfying $\|F(x^T)\|_{x^T} \leq \epsilon$ is bounded by 
\begin{equation*}
O\left(\kappa^2\log\left(\frac{\ell d_\MCal(x^0,x^\star)}{\epsilon}\right)\right). 
\end{equation*}
For each iteration where we do not repeat, the number of repeated iteration is bounded by $\lceil \log(\frac{\mu}{\hat{\mu}_0})\rceil+\lceil\log(\frac{\hat{\ell}_0}{\ell})\rceil$. Thus, the number of Riemannian gradient evaluations is bounded by
\begin{equation*}
O\left(\kappa^2\left(\log\left(\frac{\mu}{\hat{\mu}_0}\right)+\log\left(\frac{\widehat{\ell}_0}{\ell}\right)\right)\log\left(\frac{\ell d_\MCal(x^0,x^\star)}{\epsilon}\right)\right). 
\end{equation*}
This completes the proof. 
\end{proof}
\begin{remark} 
Although Algorithm~\ref{alg:FARGD} is fully adaptive, Theorem~\ref{Thm:FARGD} highlights that its convergence bound depends on $\mu$ and $\ell$. A key takeaway here is that, despite some extra terms, the fully adaptive method can retain a linear, geometry-agnostic, and last-iterate convergence rate, demonstrating its independence of any sectional curvature parameters. 
\end{remark}

%% file: sec/exp.tex
\section{Experiments}\label{sec:exp}
We present the experiments on the task of robust principal component analysis (RPCA) for symmetric positive definite (SPD) matrices. The implementation of all methods are based on the \textsc{manopt} package~\citep{Boumal-2014-Manopt} (released under GNU GPL-3 license). All the experiments were implemented in MATLAB R2024b (released under proprietary commercial software license) on a workstation with a 2.6 GHz Intel Core i7 and 16GB of memory. 
\begin{figure*}[!t]
\centering
\caption{\footnotesize{Comparison of last and average iterates for RGD with $d \in \{25, 50, 100\}$ when $(n, \alpha)=(40, 1.0)$ (above) and $(n, \alpha)=(40, 2.0)$ (bottom). The horizontal and vertical axes represent the number of data passes and the norm of Riemannian gradient.}} \label{fig:exp-RGD}
\includegraphics[width=0.3\textwidth]{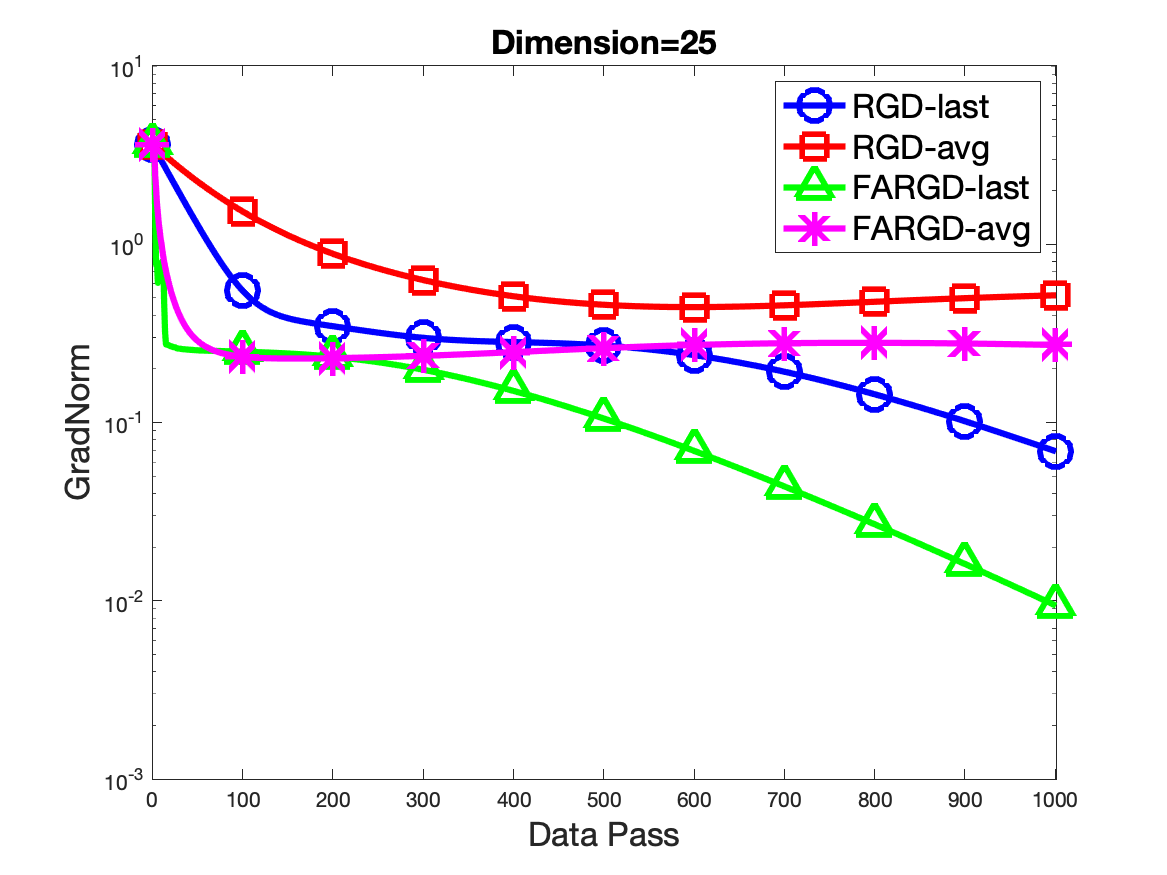}
\includegraphics[width=0.3\textwidth]{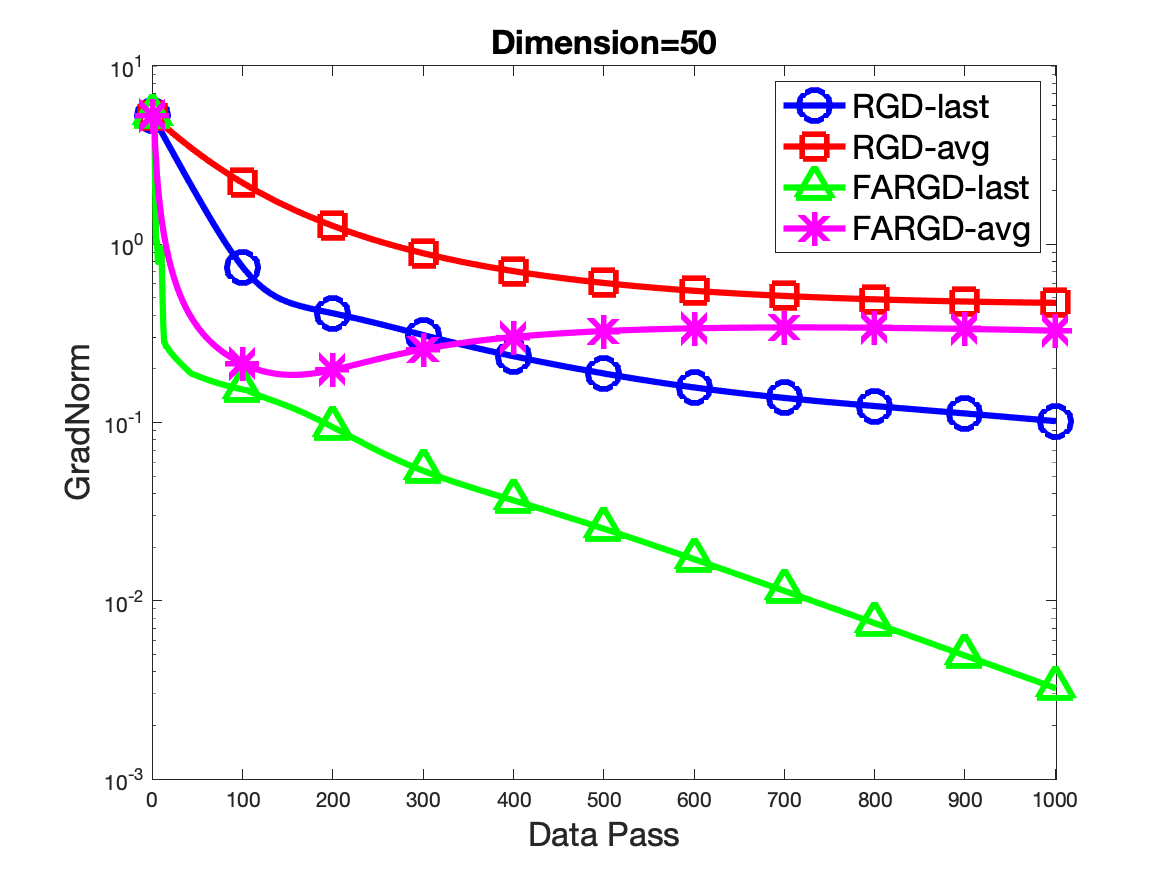}
\includegraphics[width=0.3\textwidth]{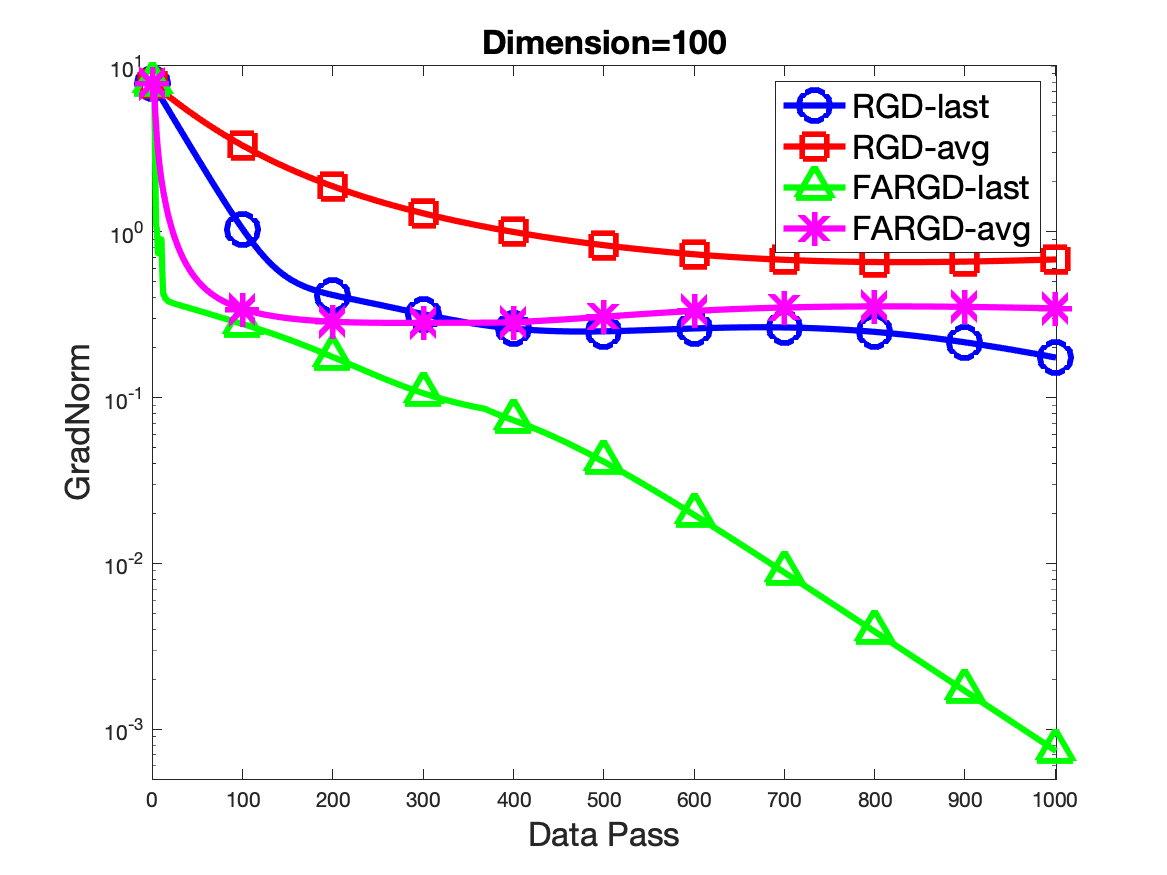} \\
\includegraphics[width=0.3\textwidth]{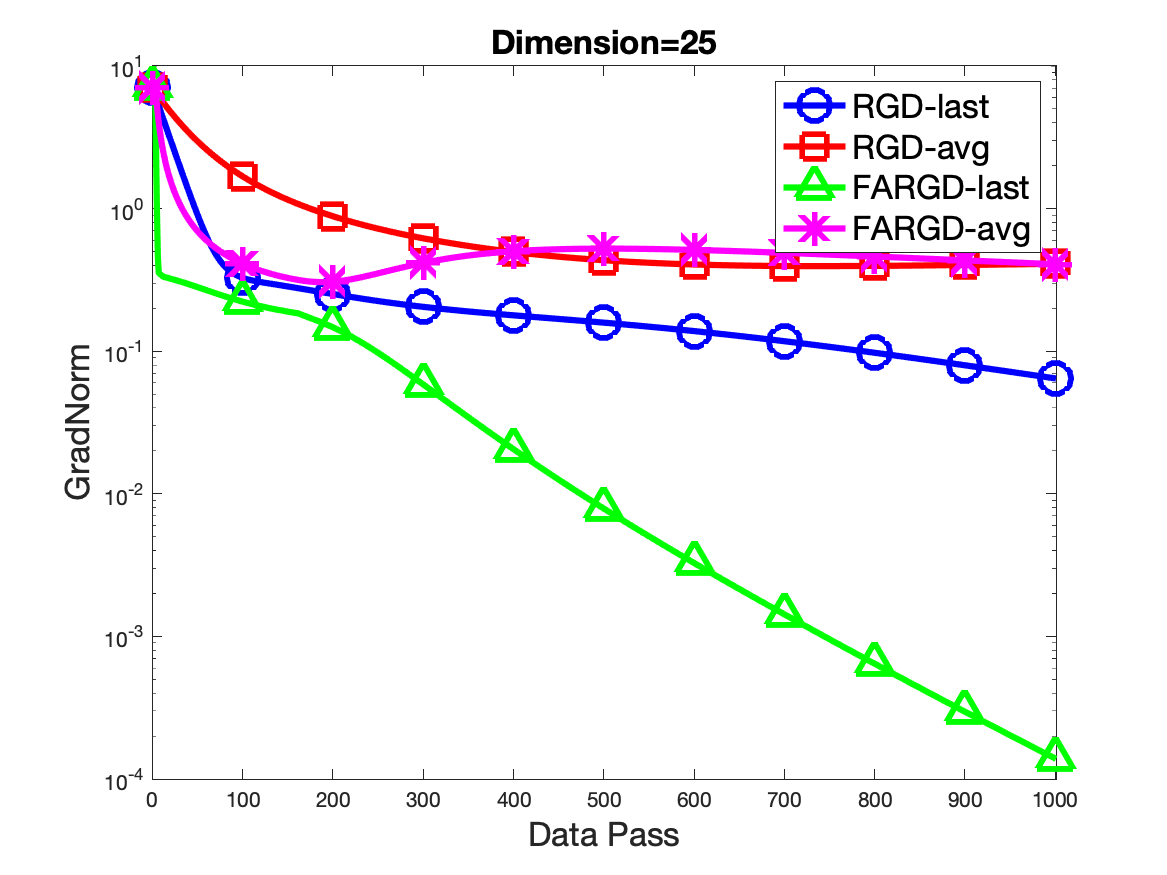}
\includegraphics[width=0.3\textwidth]{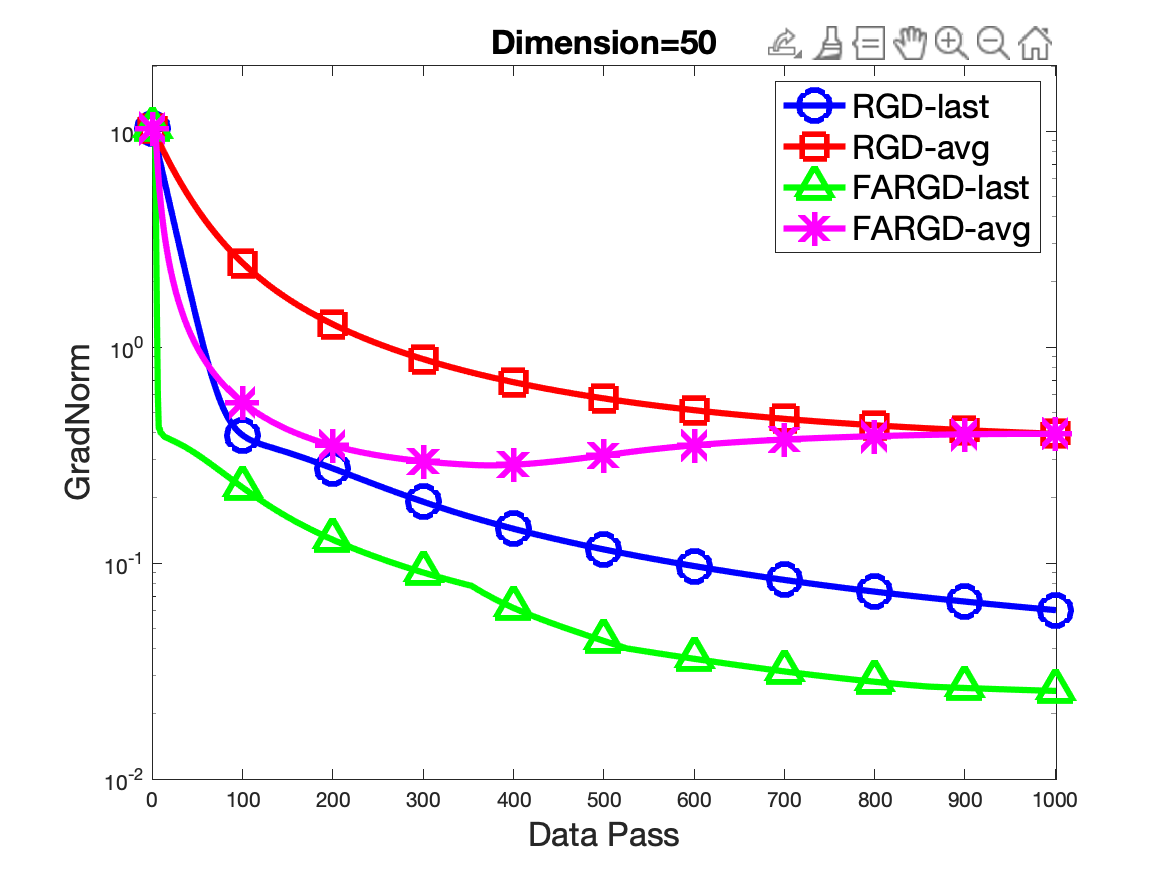}
\includegraphics[width=0.3\textwidth]{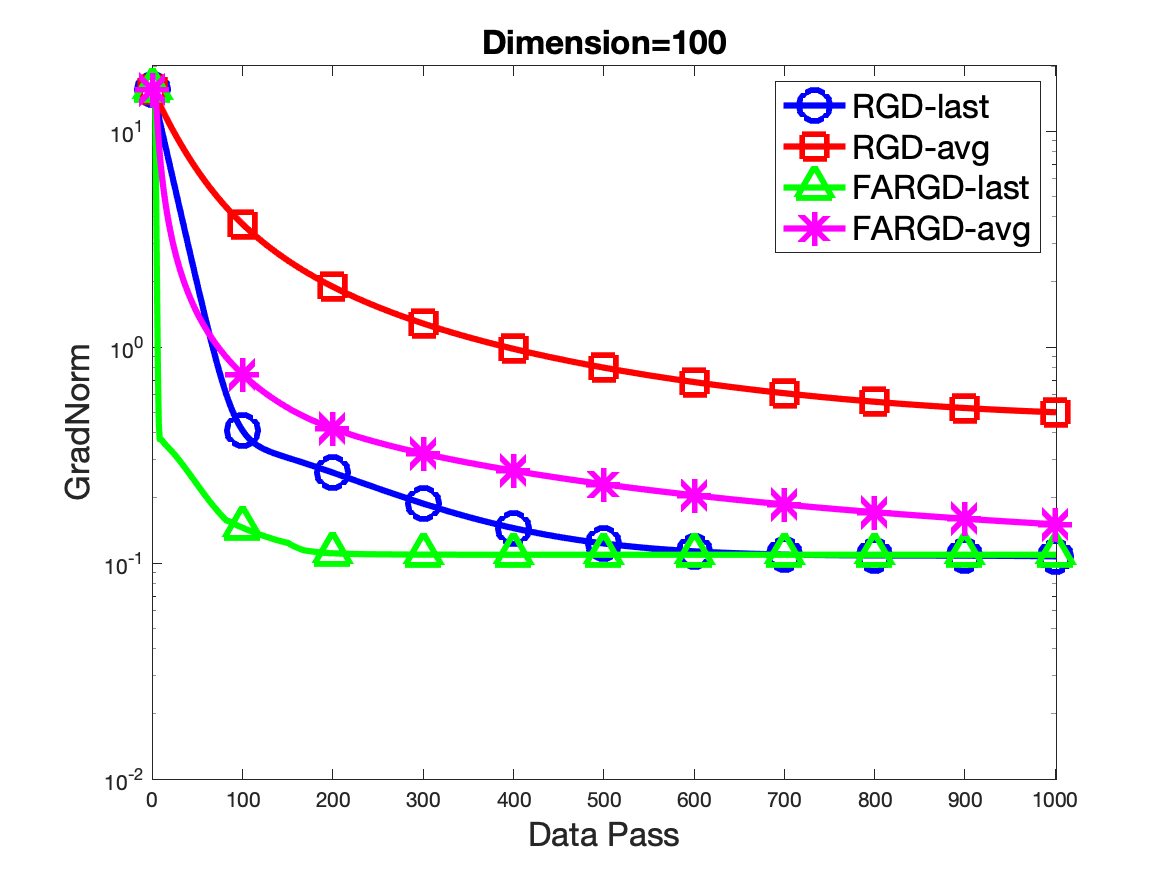}
\vspace*{-1em}
\end{figure*}
\paragraph{Experimental setup} The problem of RPCA~\citep{Candes-2011-Robust, Harandi-2017-Dimensionality} can be formulated as the Riemannian min-max optimization problem with an symmetric positive definite (SPD) manifold and a sphere manifold. Formally, we have
\begin{equation}\label{prob:RPCA}
\max_{M \in \mathcal{M}_{\textnormal{PSD}}^d} \min_{x \in \mathcal{S}^d}\left\{ -x^{\top}Mx - \tfrac{\alpha}{n}\sum_{i=1}^{n}d(M, M_i)\right\}. 
\end{equation}
In this formulation, $\alpha > 0$ is the penalty parameter, and $\{M_i\}_{1 \leq i \leq n}$ denotes a sequence of given matrices, each lying on the SPD manifold $\MCal_{\textnormal{PSD}}^d = \{M \in \br^{d \times d}: M = M^\top, M \succ 0\}$. We consider the unit sphere manifold $\SCal^d = \{x \in \br^d: \|x\| = 1\}$. The metric $d(\cdot, \cdot): \MCal_{\textnormal{PSD}}^d \times \MCal_{\textnormal{PSD}}^d \mapsto \br$ denotes the Riemannian distance induced on the SPD manifold.~\citet{Zhang-2023-Sion} has shown that, while the RPCA formulation is nonconvex-nonconcave in Euclidean space, it is locally geodesically strongly monotone. Importantly, both the SPD and sphere manifolds are \textit{complete} and \textit{boundary-free}, guaranteeing that Assumptions~\ref{Assumption:main-1} or~\ref{Assumption:main-2} hold. We highlight this example because RPCA not only serves as a standard benchmark in machine learning, but provides a clean setting to illustrate the linear convergence behavior observed in recent works~\citep{Zhang-2023-Sion, Jordan-2022-First}.

\begin{figure*}[!t]
\centering
\caption{\footnotesize{Comparison of RGD, SRGD and FARGD (last iterate) with $d \in \{25, 50\}$ when $(n, \alpha)=(200, 1.0)$. The horizontal and vertical axes represent the number of data passes and the norm of Riemannian gradient.}} 
\label{fig:exp-SRGD}
\includegraphics[width=0.4\textwidth]{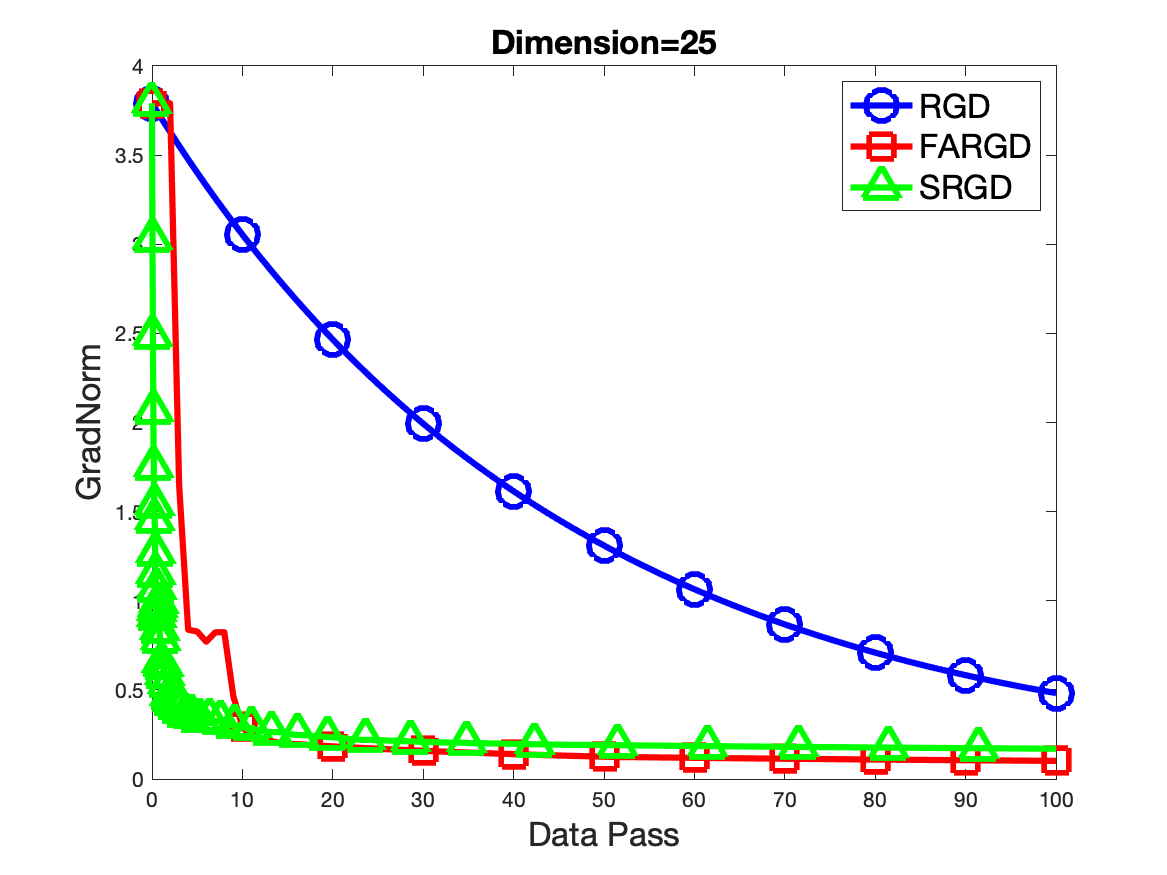}
\includegraphics[width=0.4\textwidth]{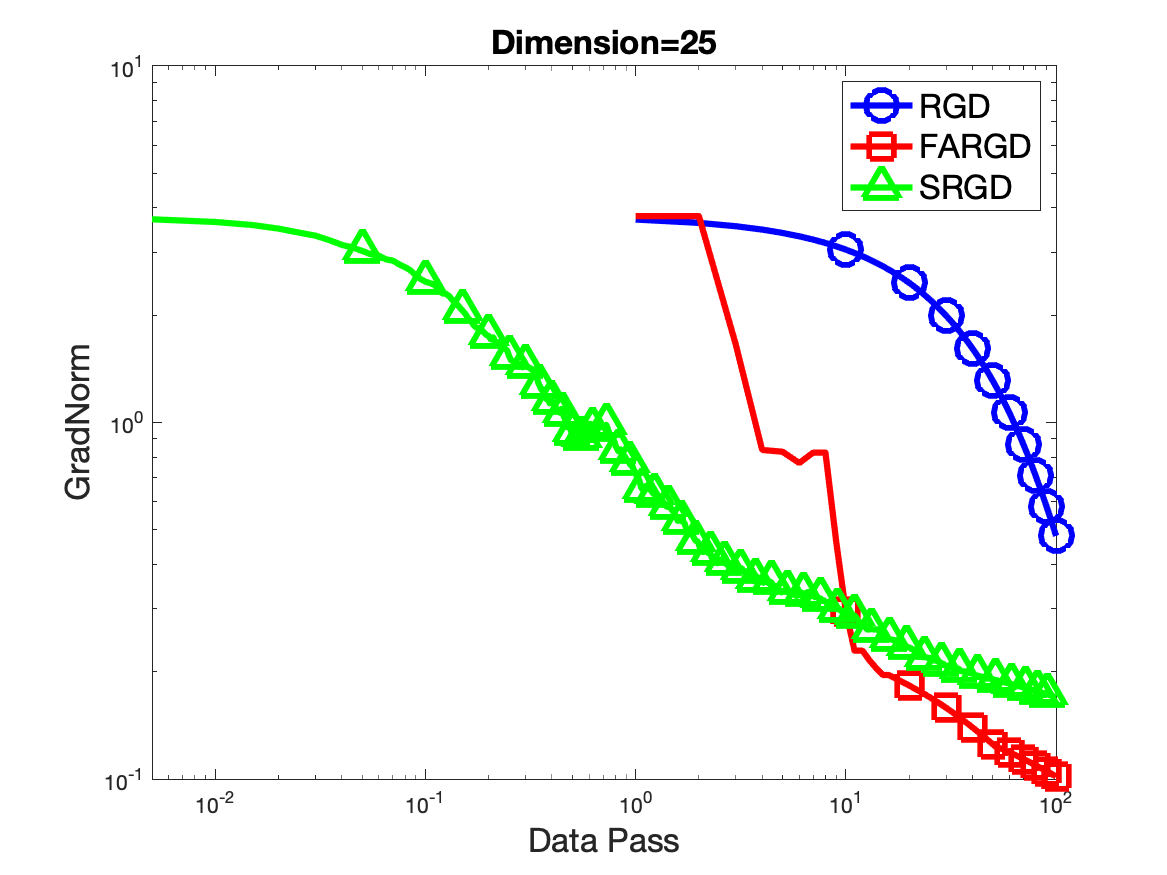} \\
\includegraphics[width=0.4\textwidth]{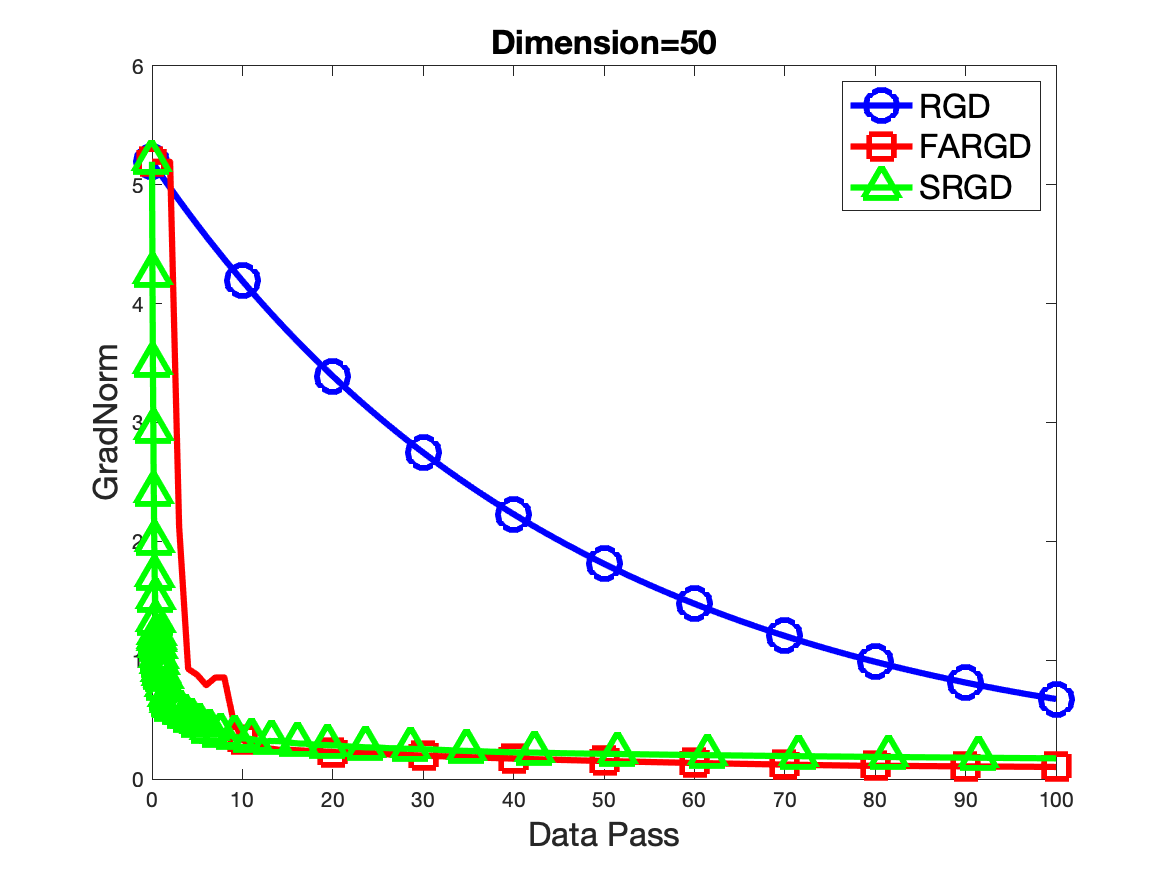}
\includegraphics[width=0.4\textwidth]{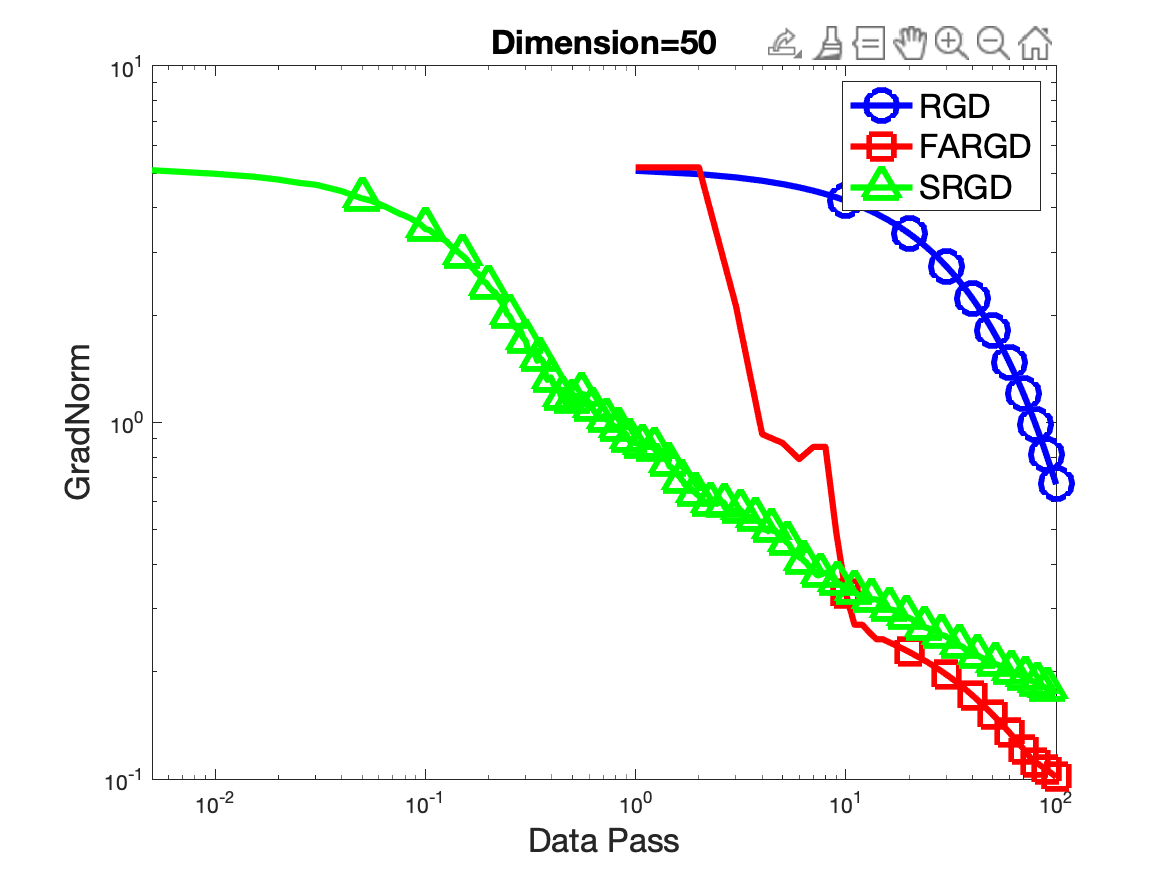}
\vspace*{-1em}
\end{figure*}

Following~\citet{Zhang-2023-Sion} and~\citet{Han-2023-Riemannian}, we generate several matrices satisfying that their eigenvalues are in the range of $[0.2,  4.5]$. In our experiment, we vary $\alpha \in \{1.0, 2.0\}$ and the problem dimension $d \in \{25, 50, 100\}$. We set the evaluation metric as the norm of Riemannian gradient and also set $\eta_t = 0.01$ (RGD and SRGD) and $(\hat{\mu}_0, \hat{\ell}_0)=(1, 1)$ (FARGD). The results with $n \in \{40, 200\}$ are summarized in Figure~\ref{fig:exp-RGD} and~\ref{fig:exp-SRGD}. 
\begin{figure*}[!t]
\centering
\caption{\footnotesize{Comparison of RGD using different stepsize choices with $d \in \{25, 50\}$ when $(n, \alpha)=(40, 1.0)$. The horizontal and vertical axes represent the number of data passes and the norm of Riemannian gradient.}} \label{fig:exp-stepsize}
\includegraphics[width=0.41\textwidth]{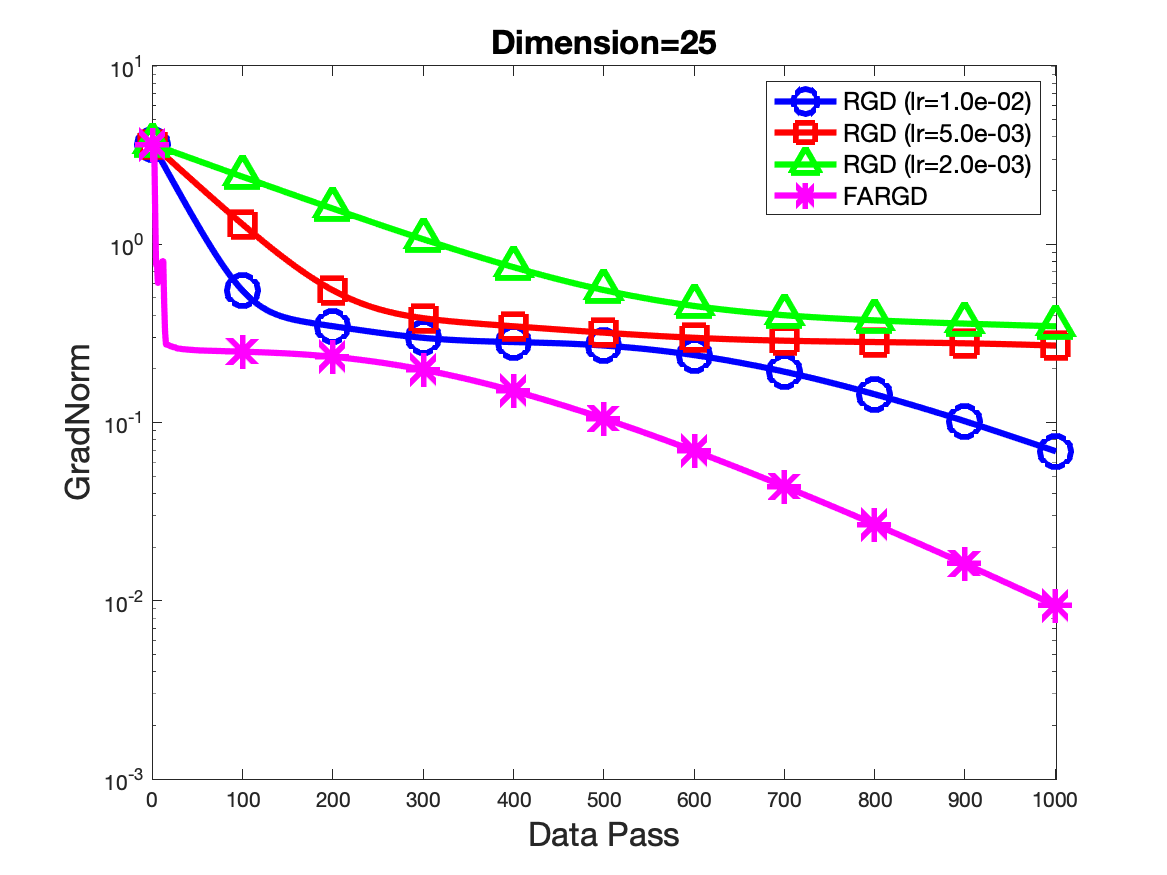}
\includegraphics[width=0.41\textwidth]{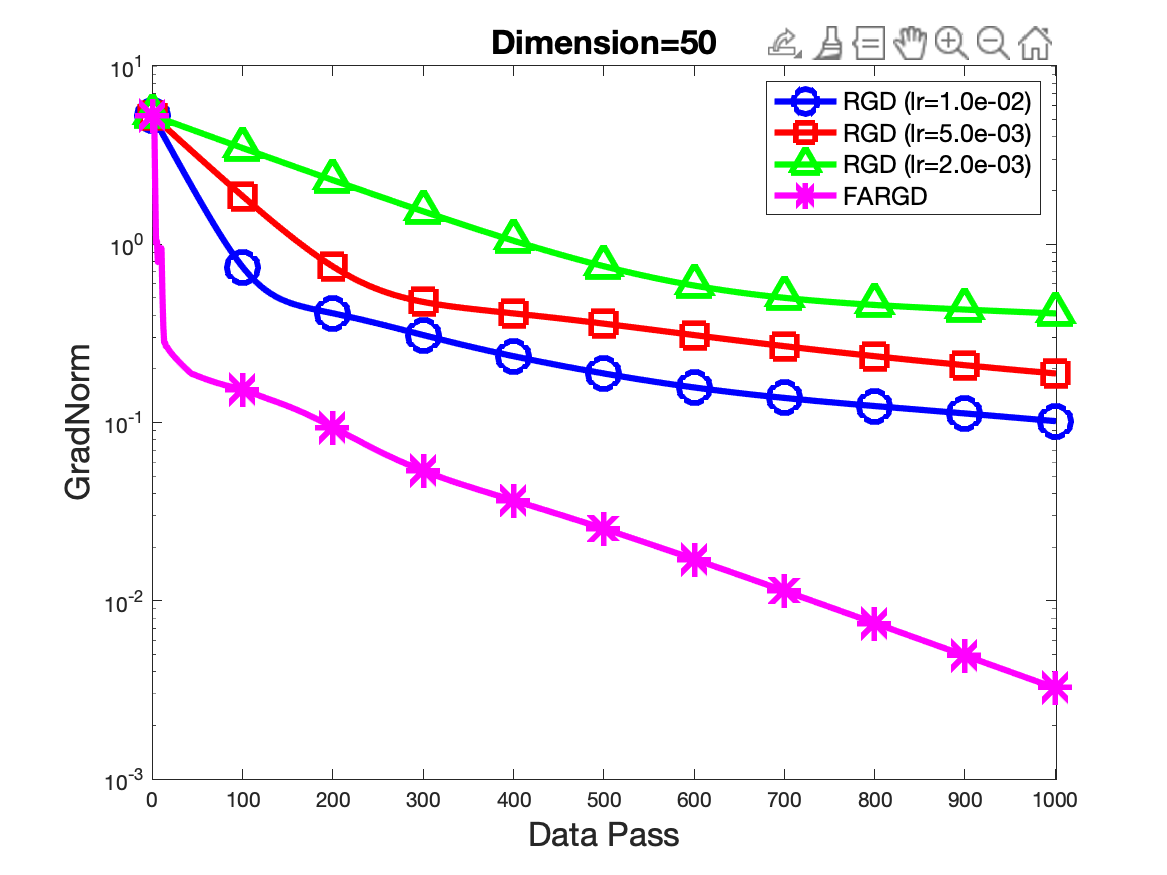}
\vspace*{-1em}
\end{figure*}
\paragraph{Experimental results} Figure~\ref{fig:exp-RGD} shows the performance of RGD and FARGD under different outputs. In particular, RGD-last/RGD-avg and FARGD-last/FARGD-avg denote the last and average iterates generated by RGD and FARGD, respectively. The results show that the last iterate of both RGD and FARGD consistently achieves a linear convergence to the optimal solution across all settings, confirming our main theoretical findings. In contrast, the average iterate converges more slowly, possibly because the RPCA problem is only locally geodesically strongly-convex-strongly-concave, and the early-stage iterates hinder the convergence speed of the average iterates. Figure~\ref{fig:exp-SRGD} compares RGD, SRGD, and FARGD all using last iterates. We observe that SRGD converges faster than RGD and FARGD in the initial phase, while all methods eventually reach the optimal solution. Finally, Figure~\ref{fig:exp-stepsize} summarizes the effect of different stepsize choices in RGD and compares them with FARGD. The results indicate that RGD is robust under various constant stepsize settings, whereas FARGD can further benefit from adaptively fitting the geometry of the problem.

%% file: sec/conclu.tex
\section{Concluding Remarks}\label{sec:conclu}
We introduced the Riemannian game framework which generalizes min-max Riemannian optimization and demonstrated that Riemannian gradient descent (RGD) can achieve a geometry-agnostic, linear and last-iterate convergence rate using a stepsize that is entirely independent of the manifold curvature. We extended these results to stochastic and fully adaptive variants of RGD. Our work opens the directions for future research. A key question is whether one can obtain best-of-both-worlds convergence in geodesically strongly monotone Riemannian games by combining geometry agnosticism with linear dependence on the condition number. Another challenge lies in establishing geometry-agnostic and last-iterate convergence for monotone Riemannian games. It also remains an open problem to design first-order algorithms that converge in structurally rich yet geodesically nonmonotone Riemannian games. Finally, our current analysis is not valid when the manifolds have boundaries. Extending geometry-agnostic convergence to such settings would require a well-defined projection mapping and constitutes an intriguing direction for future work.

%% file: sec/app.tex
\section{Additional Related Work} \label{sec:related-work}
Riemannian optimization~\citep{Absil-2009-Optimization, Boumal-2023-Introduction} is a broad and active research area that aims to develop Riemannian gradient-based algorithms exhibiting properties analogous to their well-studied Euclidean counterparts. These include deterministic methods~\citep{Ferreira-1998-Subgradient, Zhang-2016-First, Bento-2017-Iteration, Boumal-2019-Global, Criscitiello-2022-Accelerated, Kim-2022-Accelerated}, projection-free methods~\citep{Weber-2022-Projection, Weber-2022-Riemannian}, stochastic and adaptive methods~\citep{Bonnabel-2013-Stochastic, Becigneul-2019-Riemannian, Kasai-2019-Riemannian}, variance-reduced methods~\citep{Zhang-2016-Riemannian, Tripuraneni-2018-Averaging, Sato-2019-Riemannian}, and methods that escape saddle points~\citep{Criscitiello-2019-Efficiently, Sun-2019-Escaping}, among others. Within the above framework, the Riemannian gradient descent (RGD) method~\citep{Zhang-2016-First} achieves geometry-agnostic and last-iterate convergence rates.

In the unconstrained setting, if we further assume that the game is strongly monotone or that the payoff matrix in a bilinear game has all singular values bounded away from zero, gradient-based algorithms are known to achieve linear convergence rates~\citep{Daskalakis-2018-Training, Gidel-2019-Variational, Liang-2019-Interaction, Zhang-2020-Convergence, Mokhtari-2020-Unified}. By contrast, results for the constrained setting remain relatively scarce. Most existing convergence guarantees for monotone games are asymptotic~\citep{Daskalakis-2019-Last, Lei-2021-Last}. An exception is~\citet{Wei-2021-Linear}, which establishes a linear convergence rate for the optimistic gradient method in bilinear games with polytope constraints. However, their results depend on a problem-specific parameter that may approach zero, rendering the guarantee vacuous for general smooth and monotone games. More recently,~\citet{Cai-2022-Finite} derived the first last-iterate convergence rates for all smooth and monotone games, matching the known lower bounds~\citep{Golowich-2020-Last, Golowich-2020-Tight}.

From a technical viewpoint, the choice of performance measure plays a key role in pursuing last-iterate convergence. Classical gradient-based algorithms are well known to achieve a time-averaged convergence rate of $O(1/t)$ with respect to the gap function in smooth monotone games~\citep{Nemirovski-2004-Prox, Mokhtari-2020-Convergence, Kotsalis-2022-Simple}. Beyond the gap function, convergence can also be quantified using the norm of the operator in unconstrained settings~\citep{Kim-2021-Accelerated, Yoon-2021-Accelerated, Lee-2021-Fast}, or by the natural residual (and related notions) in constrained settings that satisfy additional conditions. Recently, several works~\citep{Cai-2023-Accelerated, Cai-2023-Doubly, Cai-2024-Accelerated} have demonstrated how to achieve the optimal $O(1/t)$ last-iterate convergence rate in constrained settings without relying on the cocoercivity assumption. These results have further been extended to a broader class of structured non-monotone problems, referred to as the \textit{comonotone} setting~\citep{Cai-2023-Accelerated}.

\section{Additional Proofs}
\paragraph{Proof of Lemma~\ref{Lemma:equilibrium}} We use a proof by contradiction. Indeed, suppose that $\|F(x^\star)\|_{x^\star} \neq 0$. Then, there exists player $i$ such that $\|\grad_{x_i} L_i(x^\star)\|_{x_i^\star} \neq 0$. This motivates us to consider the deviation for player $i$ with $\eta > 0$ as follows, 
\begin{equation*}
x_i = \Exp_{x_i^\star}(-\eta \cdot \grad_{x_i} L_i(x^\star)). 
\end{equation*}
Since $L_i$ is $\ell$-smooth, we have 
\begin{eqnarray*}
L_i(x_i,x_{-i}^\star) & \leq & L_i(x^\star) + \langle \grad_{x_i} L_i(x^\star), -\eta \cdot \grad_{x_i} L_i(x^\star)\rangle_{x_i^\star} + \tfrac{\ell \eta^2}{2} \|\grad_{x_i} L_i(x^\star)\|_{x_i^\star}^2 \\
& = & L_i(x^\star) - \eta\|\grad_{x_i} L_i(x^\star)\|_{x_i^\star}^2 + \tfrac{\ell \eta^2}{2}\|\grad_{x_i} L_i(x^\star)\|_{x_i^\star}^2. 
\end{eqnarray*}
This together with a sufficiently small $\eta$ yields $L_i(x_i,x_{-i}^\star) < L_i(x^\star)$, thus contradicting that the joint action profile $x^\star$ is a Nash equilibrium.

\paragraph{Proof of Lemma~\ref{Lemma:minimax}} We have $F(x_1,x_2)=(\grad_{x_1} f(x_1,x_2), -\grad_{x_2} f(x_1,x_2))$ and show that $F$ is geodesically $\mu$-strongly monotone. 

Note that $f(x_1,x_2)$ is geodesically $\mu$-strongly convex in $x_1$ and geodesically $\mu$-strongly concave in $x_2$. Then, for any two geodesics $\gamma_1: [0,1] \mapsto \MCal_1$ and $\gamma_2: [0,1]\mapsto \MCal_2$ satisfying $(\gamma_1(0),\gamma_2(0)) = (z_1, z_2)$, $(\gamma_1(1),\gamma_2(1)) = (y_1, y_2)$ and $(\dot{\gamma}_1(0),\dot{\gamma}_2(0)) = (v_1, v_2)$, we have
\begin{eqnarray*}
f(y_1,z_2) & \geq & f(z_1,z_2) + \langle \grad_{x_1} f(z_1,z_2), v_1 \rangle_{z_1} + \tfrac{\mu}{2} \| v_1 \|_{z_1}^2, \\
-f(z_1,y_2) & \geq & -f(z_1,z_2) - \langle \grad_{x_2} f(z_1,z_2), v_2 \rangle_{z_2} + \tfrac{\mu}{2} \|v_2\|_{z_2}^2.
\end{eqnarray*}
Similarly, we have:
\begin{eqnarray*}
f(z_1,y_2) & \geq & f(y_1,y_2) - \langle \grad_{x_1} f(y_1,y_2), \Gamma_{z_1}^{y_1} v_1 \rangle_{y_1}+ \tfrac{\mu}{2} \|v_1\|_{z_1}^2, \\
-f(y_1,z_2) & \geq & -f(y_1,y_2) + \langle \grad_{x_2} f(y_1,y_2), \Gamma_{z_2}^{y_2} v_2  \rangle_{y_2} + \tfrac{\mu}{2} \|v_2\|_{z_2}^2.
\end{eqnarray*}
Adding all four inequalities yields 
\begin{eqnarray*}
\lefteqn{-\langle \grad_{x_1} f(z_1,z_2), v_1 \rangle_{z_1} + \langle \grad_{x_1} f(y_1,y_2), \Gamma_{z_1}^{y_1} v_1 \rangle_{y_1}} \\ 
& & + \langle \grad_{x_2} f(z_1,z_2), v_2 \rangle_{z_2} - \langle \grad_{x_2} f(y_1,y_2), \Gamma_{z_2}^{y_2} v_2 \rangle_{y_2} \ \geq \ \mu(\|v_1\|_{z_1}^2 + \|v_2\|_{z_2}^2).
\end{eqnarray*}
Taking the parallel transport with respect to $\gamma_1$ and $\gamma_2$ yields
\begin{eqnarray*}
\langle \grad_{x_1} f(y_1,y_2),\Gamma_{z_1}^{y_1} v_1 \rangle_{y_1} & = & \langle \Gamma_{y_1}^{z_1}\grad_{x_1} f(y_1,y_2),\Gamma_{y_1}^{z_1} \Gamma_{z_1}^{y_1} v_1 \rangle_{z_1} \ = \ \langle \Gamma_{y_1}^{z_1}\grad_{x_1} f(y_1,y_2),v_1 \rangle_{z_1}, \\ 
\langle \grad_{x_2} f(y_1,y_2),\Gamma_{z_2}^{y_2} v_2 \rangle_{y_2} & = & \langle \Gamma_{y_2}^{z_2}\grad_{x_2} f(y_1,y_2),\Gamma_{y_2}^{z_2} \Gamma_{z_2}^{y_2} v_2 \rangle_{z_2} \ = \ \langle \Gamma_{y_2}^{z_2}\grad_{x_2} f(y_1,y_2),v_2 \rangle_{z_2}. 
\end{eqnarray*}
Putting these pieces together yields 
\begin{eqnarray*}
\lefteqn{\langle \Gamma_{y_1}^{z_1} \grad_{x_1} f(y_1,y_2) - \grad_{x_1} f(z_1,z_2), v_1 \rangle_{z_1}} \\ 
& & - \langle \Gamma_{y_2}^{z_2} \grad_{x_2} f(y_1,y_2) - \grad_{x_2} f(z_1,z_2), v_2 \rangle_{z_2} \ \geq \ \mu(\|v_1\|_{z_1}^2+\|v_2\|_{z_2}^2).
\end{eqnarray*}
We consider the geodesic $\gamma(t) = (\gamma_1(t),\gamma_2(t)):[0,1] \mapsto \MCal_1 \times \MCal_2$ and have that $\gamma(0) = (z_1, z_2)$, $\gamma(1) = (y_1, y_2)$ and $\dot{\gamma}(0) = (v_1, v_2)$. Thus, we have 
\begin{equation*}
\langle \Gamma_{(y_1,y_2)}^{(z_1,z_2)}F(y_1,y_2)-F(y_1,y_2),(v_1, v_2) \rangle_{(z_1,z_2)} \geq \mu (\|v_1\|_{z_1}^2+\|v_2\|_{z_2}^2). 
\end{equation*}
This completes the proof. 

\paragraph{Proof of Lemma~\ref{Lemma:relationship}} Using the Cauchy-Schwarz inequality, we have
\begin{equation}\label{inequality:main-first}
\textnormal{gap}(x; \sqrt{N}D) = \max_{v \in T_x \MCal: \|v\|_x \leq \sqrt{N}D} \langle F(x), v \rangle_x \leq \sqrt{N}D\| F(x)\|_x. 
\end{equation}
It suffices to show that $\textnormal{Tgap}(x; D) \leq \textnormal{gap}(x; \sqrt{N}D)$. By definition, we have 
\begin{equation}\label{inequality:main-second}
\textnormal{Tgap}(x; D) = \sum_{i \in \NCal} \left(L_i(x) - \min_{z_i \in \BB_{\MCal_i}(x_i; D)} L_i(z_i, x_{-i}) \right) = \sum_{i \in \NCal} \left(L_i(x) - L_i(x_i^\star, x_{-i})\right), 
\end{equation}
where $x_i^\star = \argmin_{z_i \in \BB_{\MCal_i}(x_i,D)} L(z_i,x_{-i})$.

Fixing $i \in \NCal$, we let $\gamma_i(t): [0,1] \mapsto \MCal_i$ be a distance-minimizing constant-speed geodesic with $\gamma_i(0) = x_i$ and $\gamma_i(1) = x_i^\star$. Then, the crucial step is to prove
\begin{equation}\label{inequality:convexity}
L_i(x_i^\star, x_{-i}) \geq L_i(x) + \langle \grad_{x_i} L_i (x) ,\dot{\gamma}_i(0) \rangle_{x_i}. 
\end{equation}
Indeed, we consider the one-player geodesic $\gamma(t) = (\gamma_i(t),x_{-i}):[0,1] \mapsto \MCal$. Then, we have 
\begin{equation*}
\tfrac{d}{dt} L_i(\gamma(t)) = \tfrac{d}{dt} L_i(\gamma_i(t),x_{-i}) = \langle \grad_{x_i} L_i(\gamma(t)), \dot{\gamma}_i(t) \rangle_{\gamma_i(t)} = \langle F(\gamma(t)), \dot{\gamma}(t) \rangle_{\gamma(t)}. 
\end{equation*}
We define $\phi(t) = \langle F(\gamma(t)), \dot{\gamma}(t) \rangle_{\gamma(t)}$ and show that this function is non-decreasing on $[0, 1]$. For $0 \leq t_1 < t_2 \leq 1$, the geodescially monotonicity of the game guarantees
\begin{equation*}
(t_2 - t_1)\langle \Gamma_{\gamma(t_2)}^{\gamma(t_1)} F(\gamma(t_2)) - F(\gamma(t_1)), \dot{\gamma}(t_1) \rangle_{\gamma(t_1)} \geq 0. 
\end{equation*}
Since $t_2 > t_1$, we have 
\begin{equation*}
\langle \Gamma_{\gamma(t_2)}^{\gamma(t_1)} F(\gamma(t_2)) - F(\gamma(t_1)), \dot{\gamma}(t_1) \rangle_{\gamma(t_1)} \geq 0. 
\end{equation*}
We also have 
\begin{equation*}
\langle \Gamma_{\gamma(t_2)}^{\gamma(t_1)} F(\gamma(t_2)), \dot{\gamma}(t_1) \rangle_{\gamma(t_1)} = \langle F(\gamma(t_2)),  \Gamma_{\gamma(t_1)}^{\gamma(t_2)} \dot{\gamma}(t_1) \rangle_{\gamma(t_2)} = \langle F(\gamma(t_2)), \dot{\gamma}(t_2) \rangle_{\gamma(t_2)}. 
\end{equation*}
Putting these pieces together yields that $\phi(t)$ is non-decreasing on $[0, 1]$ and hence the function $L_i(\gamma(t))$ is convex in $t$. This yields Eq.~\eqref{inequality:main-second}. 

Since $\gamma_i$ is a distance-minimizing constant-speed geodesic with $\gamma_i(0) = x_i$ and $\gamma_i(1) = x_i^\star$, we have $d_{\MCal_i}(x_i, x_i^\star) = \|\dot{\gamma}_i(0)\|$ and $\|\dot{\gamma}_i(0)\| \leq D$. Plugging Eq.~\eqref{inequality:convexity} into Eq.~\eqref{inequality:main-second} yields 
\begin{equation*}
\textnormal{Tgap}(x; D) \leq \sum_{i \in \NCal} \left(\langle \grad_{x_i} L_i (x),-\dot{\gamma}_i(0) \rangle_{x_i}\right) \leq \max_{v \in T_x \MCal: \|v\|_x \leq \sqrt{N}D} \langle F(x) , v \rangle_x = \textnormal{gap}(x; \sqrt{N}D). 
\end{equation*}
This completes the proof. 